\documentclass[12pt,a4paper]{amsart}
\usepackage{amscd}
\usepackage[all]{xy}
\usepackage{latexsym,amssymb}

\newtheorem{theorem}{Theorem}[section]
\newtheorem{prop}[theorem]{Proposition} 
\newtheorem{lemma}[theorem]{Lemma} 
\newtheorem{conjecture}[theorem]{Conjecture}
\newtheorem{corol}[theorem]{Corollary}
\theoremstyle{definition}
\newtheorem{defin}[theorem]{Definition} 
\newtheorem{exam}[theorem]{Example}
\newtheorem{defnot}[theorem]{Definitions and Notations} 
\theoremstyle{remark}
\newtheorem{remk}[theorem]{Remark}
\newtheorem{erem}{Remark}

\def\La{\Lambda} \def\Ga{\Gamma}       
\def\Om{\Omega}       \def\th{\theta}  \def\al{\alpha}       \def\de{\delta}
\def\be{\beta}        \def\eps{\varepsilon}  \def\ga{\gamma} \def\la{\lambda}   
\def\vi{\varphi}        \def\si{\sigma}       \def\om{\omega}
\def\ze{\zeta}  \def\gM{\mathfrak m}  \def\gP{\mathfrak p}  \def\gQ{\mathfrak q}
\def\fD{\mathbf d}  \def\fE{\mathbf e}  \def\fK{\mathbf k}  \def\bB{\mathbf B}
\def\fP{\mathbf p} \def\fB{\mathbf b}    \def\dM{\mathfrak M}   \def\dT{\mathfrak T}   \def\dN{\mathfrak N}
\def\DMO{\DeclareMathOperator}	\def\={\setminus}  \def\+{\oplus}  \def\8{\infty}
\def\iso{\simeq}   \def\mN{\mathbb N} \def\mZ{\mathbb Z}   \def\kT{\mathcal T}
\def\gnr#1{\langle\,#1\,\rangle}
\def\set#1{\left\{\,#1\,\right\}}
\def\setsuch#1#2{\left\{\,#1\,|\,#2\,\right\}}
\def\rA{\mathrm A}  \def\rD{\mathrm D} \def\rE{\mathrm E} \def\rT{\mathrm T} \def\rH{\mathrm H}
\def\kO{\mathcal O} \def\kG{\mathcal G} \def\kM{\mathcal M} \def\kF{\mathcal F}
\def\kN{\mathcal N} \def\kV{\mathcal V}   \def\kE{\mathcal E} \def\kJ{\mathcal J}
 \def\kH{\mathcal H} \def\kP{\mathcal P} \def\kL{\mathcal L}   \def\mP{\mathbb P}
\DMO{\supseteqc}{Spec}	\def\tX{\widetilde X}   
\DMO{\CM}{CM}		\def\Cm{Cohen--Mac\-au\-lay module}
\def\iff{if and only if }  \def\oc{one-to-one correspondence}
\DMO{\rk}{rk}	\def\larr{\longrightarrow}
\def\pX{\breve{X}}	\DMO{\VB}{VB}
\DMO{\Qc}{Qcoh}		\DMO{\coh}{Coh}
\def\khom{\mathop{\kH\!\mathit{om}}\nolimits}
\DMO{\sing}{Sing}	\DMO\he{ht}   \def\iarr{\stackrel\sim\rightarrow}
\DMO\cod{codim}		\def\cd{\cdot}
\def\tE{\widetilde E}	\def\tO{\widetilde\kO}
\def\tG{\widetilde\kG}	
\def\oO{\mathbf 0}	\def\Com{Cohen--Macaulay}
\DMO{\pr}{p}			
\def\DN{Definitions and Notations~\ref{defnot1} }
\DMO{\Kd}{Kr.dim}	\def\dch{^{\vee\vee}}
\DMO{\supp}{Supp}	\DMO{\card}{card}
\DMO{\pic}{Pic}		
   	\def\bt{\boxtimes}    \def\tN{\tilde N}
\DMO{\mo}{mod}		\DMO{\rest}{rest}
\def\AR{Aus\-lan\-der--Rei\-ten sequence}
\def\AQ{Aus\-lan\-der--Rei\-ten quiver}
\def\AT{Aus\-lan\-der--Reiten translate}
\def\toto{\leftrightarrows}	\def\hh#1{^{[#1]}}	
	\def\tA{\widetilde A}
\DMO{\syz}{syz}	\def\gdim{\mathop\mathrm{gl.dim}\nolimits}
\DMO{\ind}{ind}	\def\dvr{discrete valuation ring}
\DMO{\Hom}{Hom}    \DMO{\ext}{Ext}  \DMO{\chr}{char}  \DMO{\im}{Im}
\def\lst#1#2{#1_1,#1_2,\dots,#1_{#2}}  \def\lsto#1#2{#1_0,#1_1,\dots,#1_{#2}}
\def\row#1#2{(#1_1,#1_2,\dots,#1_{#2})}

\title[Cohen--Macaulay modules]{On Cohen--Macaulay modules over surface
singularities} 
\author[Drozd]{Yuriy A. Drozd*}
\thanks{*Supported by the CRDF Award UM2-2094} 
\address{Max-Plank-Institut f\"ur Mathematik and  Kyiv Taras Shevchenko
University}  
\email{yuriy@drozd.org}
\author[Greuel]{Gert-Martin Greuel$^\dag$}
\thanks{$^\dag$Supported by the DFG--Schwerpunkt ``Globale Methoden in der komlexen Geometrie''}
\author[Kashuba]{Irina Kashuba}
\address{Universit\"at Kaiserslautern}
\email{greuel@mathematik.uni-kl.de}
\email{kashuba@ime.usp.br}

\begin{document}
\maketitle 

\tableofcontents

\section*{Introduction}
 During the study of \Cm s on curve singularities
 (cf.~\cite{jac,drr,grk,dg1}) it was proved that these
 singularities split into three classes: 
\begin{itemize}
 \item
  \emph{\Com{} finite}, having only finitely many indecomposable \Cm s;
 \item
  \emph{\Com{} tame}, i.e., such that, for each fixed $\,r\,$,
  the indecomposable \Cm s of rank $\,r\,$ form a finite set of
  1-parameter families;
 \item
  \emph{\Com{} wild}, that can be characterized in two ways:
 \begin{itemize}
  \item
  ``geometrically'' as those having $\,n$-parameter families of
  non-isomorphic indecomposable \Cm s for arbitrary large $\,n\,$,
  \item
  ``algebraically'' as such that for every finitely generated algebra
  $\,\La\,$ there is an exact functor from the category of finite
  dimensional $\,\La$-modules to the category of \Cm s over this
  singularity, which maps non-isomorphic modules to non-isomorphic ones
  and indecomposable to indecomposable.
 \end{itemize}
 The latter property shows that the study of modules
 in the wild case is extremely complicated and needs an essentially
  new and highly non-trivial approach. 
\end{itemize}

 Moreover, it turned out that the above ``trichotomy'' is closely related
 to the position of a curve singularity in the well-known
 Arnold's list of singularities having good deformation
 properties \cite{arn}. Namely:
\begin{itemize}
 \item
  a singularity is \Com{} finite \iff it dominates one of the simple
  plane curves singularities $\,\rA_n,\rD_n,\rE_6,$ $\rE_7,\rE_8\,$; 
 \item
  a singularity is \Com{} tame \iff it dominates one of the ``serial''
  unimodal singularities $\,\rT_{pq}\,$.
\end{itemize}\noindent
 (Recall that a singularity $\,(X,x)\,$ \emph{dominates} $\,(Y,y)\,$ if
 there is a birational surjection $\,(X,x)\to(Y,y)\,$.)

 In \cite{esn,aus} it was proved that a \emph{normal surface
 singularity} (in characteristic $0$) is \Com{} finite \iff it is a
 \emph{quotient singularity}. In \cite{kah} it was shown that a
 \emph{simple elliptic} singularity is always \Com{} tame. In both cases
 a complete description of \Cm s over such singularities was
 obtained. Moreover, in \cite{kah} Kahn elaborated a general method
 relating \Cm s over normal surface singularities to vector bundles
 over some projective curves (usually singular and even non-reduced).

 In this paper we use the Kahn correspondence and the results of
 \cite{dg2} to prove the following main theorems:
\begin{itemize}
 \item
  every cusp singularity \cite{hir,kar} is \Com{} tame;
 \item
  a minimally elliptic singularity \cite{lau} that is neither simple
  elliptic nor a cusp one is \Com{} wild.
\end{itemize}
 As a corollary we get that any \emph{log-canonical} surface
 singularity \cite{kaw} is \Com{} tame (or finite, as quotient
 singularities are also log-canonical). There is some evidence that any
 other normal surface singularity is \Com{} wild. We also give a
 description of \Cm s over cusp singularities and use it to provide
 a description of \Cm s over curve singularities of type $\,\rT_{pq}\,$ (there
 was no such description in \cite{dg1}, the proof of their tameness was indirect).

\newpage
\section{Generalities}
\label{sec1}

 In what follows we use the following 
\begin{defnot}
\label{defnot1}
\begin{itemize}
 \item
  A \emph{surface singularity} means a spectrum $\,X=\supseteqc A\,$, where
  $\,A\,$ is a local, complete\,%
 \footnote{
  For the analytic case, see Remark~\ref{anal} at the end of the
  paper. 
  }
  noetherian ring of Krull dimension
  $2$. We  denote by $\,\gM=\gM_A\,$ the maximal ideal of $\,A\,$ and
  by $\,\fK=A/\gM\,$ the residue field. We always suppose that the
  field $\,\fK\,$ is algebraically closed.
 \item
  Such a singularity is called \emph{normal} if the ring $\,A\,$ is
  normal, i.e., integral and integrally closed in its field of
  fractions $\,Q\,$. In what follows we only consider normal surface
  singularities. Recall that any normal surface singularity is
  \emph{isolated}, that is $\,\gM\,$ is the unique singular point of
  it. 
 \item
  A \emph{resolution} of a normal surface singularity is a birational
  proper map $\,\pi:\tX\to X\,$, where $\,\tX\,$ is smooth, which
  induces an isomorphism $\,\tX\=E\iarr X\=\set{\gM}\,$, where
  $\,E=\pi^{-1}(\gM)_{\mathrm{red}}\,$. $\,E\,$ is called the \emph{exceptional
  curve} on $\,\tX\,$ (it is indeed a projective curve over
  $\,\fK\,$). Denote by $\,E_i\,$ ($\,i=1,\dots,s\,$) the irreducible
  components of $\,E\,$. Put also $\,\pX=X\=\set{\gM}\,$.
 \item
  A resolution $\,\pi\,$ as above is called \emph{minimal} if it
  cannot be decomposed as
$$
\begin{CD}
  \tX @>\pi_1>> X' @>\pi_2>> X \,,
\end{CD}
$$
 where $\,X'\,$ is also smooth. It is known (cf., e.g., \cite{lip})
 that $\,\pi\,$ is minimal \iff neither of the components $\,E_i\,$ is
 a smooth rational curve with self-intersection index
 $\,E_i\cd E_i=-1\,$. 
 \item
  A \emph{cycle} on $\,\tX\,$ is a divisor  $\,Z=\sum_{i=1}^sk_iE_i\,$ with $\,k_i\in\mZ\,$.
 Such a cycle is called \emph{positive} if
 all $\,k_i>0\,$. We write $\,Z\le Z'\,$ if
 $\,Z'=\sum_{i=1}^sk_i'E_i\,$ with $\,k_i\le k_i'\,$ for all $\,i\,$.
 \item
  The \emph{fundamental cycle}, denoted by $\,Z_0\,$, is the minimal
 positive cycle (with respect to the above defined partial order) such that
 $\,Z_0\cd E_i\le0\,$ for all $\,i\,$.
 \item
  A positive cycle $\,Z\,$ is called a \emph{weak reduction cycle} if
 the sheaf $\,\kO_Z(-Z)\,$ is \emph{generically generated} by its
 global section (i.e. generated outside a finite set) and
 $\,\rH^1(E,\kO_Z(-Z))=0\,$. Here, as usually, we write
 $\,\kO_Z=\kO_{\tX}/\kO_{\tX}(-Z)\,$ and consider $\,Z\,$ as a closed
 subscheme of $\,\tX\,$. 
 \item
  A \emph{reduction cycle} is a weak reduction cycle $\,Z\,$ such that
 the sheaf $\,\om^\vee_Z\,$ is generically generated by global
 sections, where $\,\om_Z=\om_{\tX}\*_{\kO_{\tX}}\kO_Z(-Z)\,$ is the
 dualizing sheaf for $\,Z\,$ and $\,\kM^\vee\,$ denotes
 $\,\khom_{\kO_S}(\kM,\kO_S)\,$, where $\,\kM\,$ is a
 coherent sheaf on a scheme $\,S\,$. 
 \item
  One calls a surface singularity \emph{rational} if
 $\,\rH^1(\tX,\kO_{\tX})=0\,$ and \emph{minimally elliptic} if it is
 Gorenstein (i.e. $\,\om_X\iso\kO_X\,$) and $\,\rH^1(\tX,\kO_{\tX})\iso\fK\,$. 
 \item
  $\,X\,$ is called a \emph{simple elliptic singularity} if $\,E\,$ is a
 smooth curve of genus $1$.
 \item
  We call the exceptional curve $\,E\,$ a \emph{cyclic configuration} in the following
  cases:
  \begin{itemize}
   \item[(i)]
   $\,s=1\,$, $\,E\,$ is rational and has a unique singular point
   that is a simple node;
   \item[(ii)]
   $\,s=2\,$, $\,E_1\iso E_2\iso\mP^1\,$ and they intersect
   transversally in (exactly) $2$ points;
   \item[(iii)]
   $\,s>2\,$, $\,E_i\iso\mP^1\,$ for all $\,i\,$, $\,E_i\cd
   E_{i+1}=1\,$, $\,E_s\cd E_1=1\,$, and $\,E_i\cd E_j=0\,$
   otherwise. 
  \end{itemize}
  For a cyclic configuration we set $\,E_{s+k}=E_k\,$ for all integers $\,k\,$. 
 \item
  $\,X\,$ is called a \emph{cusp singularity} if $\,E\,$ is a cyclic
 configuration. Note that simple elliptic and cusp singularities
 are both minimally elliptic. Moreover, for all of them $\,Z_0=E\,$.
 \item 
  Denote by $\,\coh(S)\,$ the category of coherent sheaves on a scheme
 $\,S\,$ and by $\,\VB(S)\,$ the category of \emph{vector bundles},
 that is locally free coherent sheaves on $\,S\,$.
 \item
  For a surface singularity $\,X=\supseteqc A\,$ denote by $\,\CM(X)\,$, or
 by $\,\CM(A)\,$, the category of (maximal) \Cm s over the ring
 $\,A\,$. If $\,A\,$ is normal, \Cm s coincide with \emph{reflexive}
 modules, i.e. such $\,A$-modules $\,M\,$ that $\,M\dch\iso
 M\,$, where $\,M^\vee=\Hom_A(M,A)\,$. We always identify an
 $\,A$-module $\,M\,$ with its ``sheafification'' $\,\tilde M\,$, which
 is a quasi-coherent sheaf on $\,X\,$.
\end{itemize}
\end{defnot}

 Recall the main result of Kahn \cite[Theorem 1.4]{kah} concerning the
 relations between the vector bundles on a reduction cycle and the \Cm s
 over a normal surface singularity. 
\begin{theorem}
\label{kahn1}
 Let $\,\pi:\tX\to X\,$ be a resolution of a normal surface
 singularity with exceptional curve $\,E\,$ and a reduction cycle
 $\,Z\,$. Denote by $\,R_Z:\CM(X)\to\VB(Z)\,$ the functor that is the
 composition of the functors
\begin{align*}
 \CM(X)\to\VB(\tX) &:\ M \mapsto (\pi^*M)\dch\\
\intertext{and}
 \VB(\tX)\to\VB(Z) &:\ \kF\mapsto\kF\*_{\kO_{\tX}}\kO_Z\,.
\end{align*}
 This functor maps non-isomorphic objects to non-isomorphic ones, and a
 vector bundle $\,\kV\in\VB(Z)\,$ is isomorphic to $\,R_ZM\,$ for some
 $\,M\,$ \iff it is generically generated by global sections and there
 is an extension of $\,\kV\,$ to a vector bundle $\,\kV_2\,$ on
 $\,2Z\,$ such that the exact sequence
$$
 0\larr \kV(-Z)\larr \kV_2\larr \kV\larr 0 \,,
$$
 induces a monomorphism $\,\rH^0(E,\kV(Z))\to \rH^1(E,\kV)\,$.
\end{theorem}

 We denote by $\,\VB^K(Z)\,$ the full subcategory of $\,\VB(Z)\,$
 consisting of the vector bundles satisfying the latter Kahn's
 conditions. Note that obviously $\,\rk R_ZM=\rk M\,$ for any \Cm{}
 $\,M\,$. 

 For minimally elliptic singularities one can give a simpler
 description of the latter category (cf. \cite[Theorem
 2.1]{kah}).
\begin{theorem}
\label{kahn2}
 Let $\,\pi:\tX\to X\,$ be a minimal resolution of a minimal elliptic
 singularity. Then the fundamental cycle $\,Z=Z_0\,$ is a reduction
 cycle and the category $\,\VB^K(Z)\,$ consists of all vector bundles
 of the form $\,n\kO_Z\+\kG\,$, where $\,\kG\,$ satisfies the following
 conditions:
\begin{enumerate}
 \item
  $\,\kG\,$ is generically generated by global sections,
 \item
  $\,\rH^1(E,\kG)=0\,$,
 \item
  $\,\dim_\fK\rH^0(E,\kG(Z))\le n\,$.
\end{enumerate}
 In particular, indecomposable objects of $\,\VB^K(Z)\,$ are the
 following:
\begin{itemize}
 \item
  the trivial line bundle $\,\kO_Z\,$,
 \item
  $\,n\kO_Z\+\kG\,$, where $\,\kG\,$ is an indecomposable vector
  bundle satisfying the above conditions {\rm (1),(2)}, and
  $\,n=\dim_\fK\rH^0(E,\kG(Z))\,$. 
\end{itemize}
\end{theorem}

\begin{erem}
 In the paper \cite{kah}, Theorems \ref{kahn1} and \ref{kahn2} are
 proved in the ``geometric'' case, when $\,A\,$ is a $\,\fK$-algebra
 and $\,\chr\fK=0\,$, but one easily verifies that all proofs can be
 directly extended to the general situation. The results of Laufer
 \cite{lau} used in \cite{kah} are also valid in this ``abstract''
 context. Using the Grauert--Riemenschneider vanishing theorem to show that
 $\,\rH^1(E,\kO_Z(-Z))=0\,$ for the fundamental cycle of a minimally elliptic singularity, is actually superfluous.
 Indeed, in this case $\,\om_{\tX}=\kO_{\tX}(-Z)\,$ \cite[Theorem 3.4]{lau}, so $\,\om_Z\iso\kO_Z\,$.
 By the Serre's duality \cite[Corollary III.7.7]{har}, $\,\rH^1(E,\kO_Z(-Z))\,$ is dual to
 $\,\rH^0(E,\om_Z(Z))=\rH^0(E,\kO_Z(Z))\,$. The latter is zero since $\,\deg\kO_Z(Z)=Z\cdot Z<0\,$.
\end{erem}

\section{\Cm s on cusp singularities}
\label{sec2}

 In this section $\,E\,$ denotes a cyclic configuration,
 $\,\nu:\tE\to E\,$ the normalization of $\,E\,$ and $\,E_i\ (i=1,\dots,s)\,$
 the irreducible components of $\,\tE\,$ (if $\,s>1\,$, they can be identified
 with the irreducible components of $\,E\,$). We put
 $\,E_{s+k}=E_k\,$ for all integers $\,k\,$, $\,\kO=\kO_E\,$,
 $\,\tO=\kO_{\tE}\,$, $\,\kO_i=\kO_{E_i}\,$, and we always identify
 $\,\tO\,$ and $\,\kO_i\,$ with their images under $\,\nu_*\,$. Denote by
 $\,S=\set{\lst ps}\,$ the set of singular points of $\,E\,$ and choose
 the indices in such a way that $\,\nu^{-1}(p_i)\subset E_i\cup
 E_{i+1}\,$. Note that the latter preimage always consists of $2$
 points, which we denote by $\,p'_i\,,\,p''_{i+1}\,$ so that
 $\,p'_i\in E_i\,$ and $\,p''_{i+1}\in E_{i+1}\,$. Let $\,\kJ\,$ be the
 \emph{conductor} of $\,\tO\,$ in  $\,\kO\,$, i.e. the biggest
 $\,\tO$-ideal contained in $\,\kO\,$. Then $\,\kO_i/\kJ\kO_i=\fK(p'_i)\+\fK(p''_i)\,$.

 Recall the description of vector bundles on a cyclic configuration
 given in \cite{dg2}. Call an \emph{$\,s$-sequence} any sequence of
 integers  of length a multiple of $\,s\,$: $\,\fD=(d_1,\dots,d_{rs})\,$.
 Call this sequence \emph{aperiodic} if it cannot be obtained
 by repetition of a shorter $\,s$-sequence. Call an \emph{$\,s$-shift}
 of $\,\fD\,$ any sequence $\,\fD^{k}=(d_{ks+1},\dots,d_{rs},d_1,\dots,d_{ks})\,$. 

\begin{theorem}
\label{dg}
 Indecomposable vector bundles on a cyclic configuration are in \oc
 with the triples $\,\bB=(\fD,m,\la)\,$, where $\,m\,$ is a positive
 integer, $\,\la\in\fK^*=\fK\=\set{0}\,$ and $\,\fD\,$ is a aperiodic
 $\,s$-sequence defined up to an $\,s$-shift.
\end{theorem}
 We denote the vector bundle corresponding to the triple $\,\bB\,$ by
 $\,\kG(\bB)\,$. The precise construction of this sheaf is the
 following:
\begin{enumerate}
\label{ei}
 \item
  Let $\,\fD=(d_1,\dots,d_{rs})\,$. Put $\,\kG_i=m\kO_i(d_i)\,$ and 
  $\,\tG(\bB)=\bigoplus_{i=1}^{rs}\kG_i\,$. Note that
  $\,\kG_i/\kJ\kG_i\iso m\fK(p'_i)\+m\fK(p''_i)\,$. Put also
  $$\,F(\bB)=\tG(\bB)/\kJ\tG(\bB)\iso\bigoplus_{i=1}^{rs}
  \bigl(m\fK(p'_i)\+m\fK(p''_i)\bigr)\,.$$
  Since $\,F(\bB)\,$ is a sky-scraper sheaf, we identify it with the
  vector space of its global sections. 
 \item
  Choose a $\,\fK$-bases  $\,\set{\fE'_{ik},\fE''_{ik}}\,$ ($\,k=1,\dots,m\,$)
 of $\,F(\bB)\,$ in such a way
  that $\,\set{\fE'_{ik}}\,$ form a basis of $\,m\fK(p'_i)\,$ and
  $\,\set{\fE''_{ik}}\,$ form a basis of $\,m\fK(p''_i)\,$.
 \item
  Define the elements $\,\fE_{ik}\in F(\bB)\,$
  ($\,i=1,\dots,rs,\,k=1,\dots,m\,$) in the following way:
$$
 \fE_{ik}=\begin{cases}
 	   \fE_{ik}'+\fE_{i+1,k}'' &\textrm{if }\ i\ne rs\,, \\
	   \fE_{rs,k}'+\la\fE_{1k}''+\fE_{1,k-1}'' &\textrm{if }\
 	   i=rs,\, k\ne1 \,,\\
           \fE_{rs,1}'+\la\fE''_{11} &\textrm{if }\ i=rs,\,k=1\,.
 	  \end{cases}
$$
 Denote by $\,G(\bB)\,$ the subspace of $\,F(\bB)\,$ with the basis 
 $$\,\setsuch{\fE_{ik}}{i=1,\dots,rs,\,k=1,\dots,m}\,.$$
 \item
  Now $\,\kG(\bB)\,$ is defined as the preimage of $\,G(\bB)\,$ in
  $\,\tG(\bB)\,$ under the epimorphism $\,\kG(\bB)\to F(\bB)\,$. 
  Note that $\,\rk\kG(\bB)=mr\,$.
\end{enumerate}

 To use Kahn's Theorem we have to calculate the cohomologies of the
 sheaves $\,\kG(\bB)\,$. Let
 $\,H(\bB)=\tG(\bB)/\kG(\bB)=F(\bB)/G(\bB)\,$. The exact sequence
 $\,0\to\kG(\bB)\to\tG(\bB)\to H(\bB)\to0\,$ of sheaves over $\,E\,$
 gives the exact sequence of cohomology groups
$$
 0\to \rH^0(\kG(\bB))\to \rH^0(\tG(\bB))\stackrel{h(\bB)}{\larr} H(\bB) \to 
 \rH^1(\kG(\bB))\to \rH^1(\tG(\bB))\to0
$$
 (we identify $\,H(\bB)\,$ with the space of its global
 sections and omit ``$E$'' in the notation of cohomology). One knows
 that $\,\dim_\fK\rH^0(\tG(\bB))=m\sum_{i=1}^{rs}(d_i+1)^+\,$ and
 $\,\dim_\fK\rH^1(\tG(\bB))=m\sum_{i=1}^{rs}(d_i+1)^-\,$, where
$$
 k^+=\begin{cases}
  k &\textrm{if }\ k>0\,,\\ 0 &\textrm{otherwise}
 \end{cases}
\qquad\textrm{and}\qquad
 k^-=\begin{cases}
  -k &\textrm{if } k<0\,,\\ 0 &\textrm{otherwise}\,,
 \end{cases}\,,
$$
 while $\,\dim_\fK H(\bB)=mrs\,$. So we only have to calculate $\,\rk
 h(\bB)\,$. Denote by $\,f\,$ the natural homomorphism
 $\,\rH^0(\tG(\bB))\to F(\bB)\,$. Then $\,\rk h(\bB)=\dim_\fK(\im
 f+G(\bB))-\dim_\fK G(\bB)\,$.

 Identifying $\,E_i\,$ with $\,\mP^1\,$ we can suppose that
 $\,p'_i=(1:0)\,$ and $\,p''_i=(0:1)\,$. If $\,d_i>0\,$, a basis
 $\,\set{g_0,\dots,g_{d_i}}\,$ in $\,\rH^0(\kO_i(d_i))\,$ can be
 chosen in such a way that $\,g_0(p'_i)=1,\,g_1(p''_i)=1\,$, while all
 other values $\,g_j(p'_i)\,$ and $\,g_j(p''_i)\,$ equal $\,0\,$, so
 the basis elements $\,\fE'_{ik}\,$ and $\,\fE''_{ik}\,$ can be chosen
 as the images of $\,g_0\,$ and $\,g_1\,$ respectively. Hence $\,\im
 f\,$ is generated by the set $\,\setsuch{\fE'_{ik},\fE_{ik}''}{d_i>0}\cup
 \setsuch{\fE'_{ik}+\fE''_{ik}}{d_i=0}\,$.

 Call a subsequence $\,\fP=(d_{k+1},\dots,d_{k+l})\,$, where $\,0\le
 k<rs\,$ and $\,1\le l\le rs\,$, a \emph{positive part}
 of $\,\fD\,$ if $\,d_{k+j}\ge0\,$ for $\,j=1,\dots,l\,$ and
 either $\,l=rs\,$ or both $\,d_k<0\,$ and $\,d_{k+l+1}<0\,$. For such
 a positive part put $\,\th(\fP)=l\,$ if either $\,l=rs\,$ or
 $\,\fP=(0,\dots,0)\,$ and $\,\th(\fP)=l+1\,$ otherwise. Denote by
 $\,\kP\,$ the set of all positive parts of $\,\fD\,$ and put
 $\,\th(\fD)=\sum_{\fP\in\kP}\th(\fP)\,$. Put also
 $\,\de(\fD,\la)=1\,$ if $\,\fD=(0,\dots,0)\,$ and $\,\la=1\,$ and
 $\,\de(\fD,\la)=0\,$ otherwise. The preceding consideration and a
 rather simple straightforward calculations show that, for
 $\,\bB=(\fD,m,\la)\,$,
$$
 \rk H(\bB)=m\th(\fD)-\de(\fD,\la)\,,
$$
 whence the following:
\begin{theorem}
\label{cohom}
 For the vector bundle $\,\kG=\kG(\fD,m,\la)\,$,
\begin{align*}
 \dim_\fK \rH^0(E,\kG) &=
 m\left(\sum_{i=1}^rs(d_i+1)^+-\th(\fD)\right)+\de(\fD,\la) \\
 \dim_\fK \rH^1(E,\kG) &=
 m\left(\sum_{i=1}^rs(d_i+1)^-+rs-\th(\fD)\right)+\de(\fD,\la)\,.
\end{align*}
\end{theorem}

 We write $\,\fD>0\,$ if $\,d_i\ge0\,$ for all $\,i=1,\dots,rs\,$ and
 at least one of them is strictly positive. Denote by $\,\oO\,$ the
 zero sequence $\,(0,\dots,0)\,$ of length $\,s\,$ (the unique
 aperiodic zero $\,s$-sequence).
\begin{corol}
\label{cusp0}
 A vector bundle $\,\kG=\kG(\fD,m,\la)\,$ satisfies Kahn's conditions
 {\rm(1),(2)} of\ {\,\rm Theorem~\ref{kahn2}} \iff either $\,\fD>0\,$ or
 $\,\fD=\oO\,$ and $\,\la\ne1\,$.
\end{corol}

 Suppose now that the cycle configuration $\,E\,$ is actually the
 exceptional curve of a cusp singularity $\,X=\supseteqc A\,$. If
 $\,s>1\,$, set $\,b_i=-E\cdot E_i=-E_i\cdot E_i-2\,$. Note
 that the intersection indices $\,E_i\cd E_i\le -2\,$ and at least one of them is smaller than
 $\,-2\,$. If $\,s=1\,$, set $\,b_1=-E\cd E\,$. Let
 $\,\fB^r=(b_1,\dots,b_s,\dots,b_1,\dots,b_s)\,$ ($\,r\,$ times). We
 write $\,\fB\,$ for $\,\fB^1\,$. Put
 also $\,n(\kG)=\dim_\fK\rH^0(E,\kG(E))\,$. Recall that in this case
 one can choose $\,E\,$ itself as a reduction cycle.
\begin{corol}
\label{cusp1}
\begin{enumerate}
 \item Suppose that a vector bundle $\,\kG=\kG(\fD,m,\la)\,$ satisfies
  the Kahn conditions {\rm(1),(2)}. Then
$$
 n(\kG)=m\left(\sum_{i=1}^{rs}(d_i-b_i+1)^+-\th(\fD-\fB^r)\right)+
 \de(\fD-\fB^r,\la),
$$
 so $\,n(\kG)\kO\+\kG\iso R_EM(\fD,m,\la)\,$, where
 $\,M(\fD,m,\la)\,$ is a  Cohen--Macaulay $\,A$-module of rank
 $\,rs+n(\kG)\,$.
 \item
  All \Cm s $\,M(\fD,m,\la)\,$ are indecomposable and every
 indecomposable Cohen--Macaulay $\,A$-module is isomorphic either to 
 $\,A\,$ or to one of the modules $\,M(\fD,m,\la)\,$.
\end{enumerate}
\rm Note that here either $\,\fD>0\,$ or $\,\fD=\oO\,$ and
$\,\la\ne1\,$. 
\end{corol}

\section{Lifting families}
\label{sec3}

 Vector bundles on cyclic configurations form natural
 $\,1$-parametric families. We are going to show that these families
 can be lifted to $\,1$-parametric families of \Cm s on cusp
 singularities. To do it, we first prove analogues of Kahn's results
 for families of modules and vector bundles. For the sake of
 simplicity, we suppose from now on that $\,A\,$ is a
 $\,\fK$-algebra and all other algebras and schemes are
 $\,\fK$-algebras and $\,\fK$-schemes. Moreover, we only consider
 $\,\fK$-schemes $\,S\,$ such that $\,\fK(x)=\fK\,$ for each closed
 point $\,x\in S\,$ (for instance, schemes of finite type over
 $\,\fK\,$). We write $\,\*,\,\dim,\,\times\,$, etc., instead of
 $\,\*_\fK,\,\dim_{\fK},\,\times_{\supseteqc\fK}\,$, etc. Especially ``finite dimensional''
 means ``finite dimensional over $\,\fK$.''

 We need also families with non-commutative base,
 so we give the corresponding definitions.

\begin{defin}\label{fam0}
 Let $\,S\,$ be a $\,\fK$-scheme and $\,\La\,$ be a $\,\fK$-algebra (maybe non-commutative). 
\begin{itemize}
\item  A \emph{family of  $\,\kO_S$-modules based on} $\,\La\,$ is, by definition, a coherent sheaf of
 $\,\kO_S\*\La$-modules $\,\kF\,$ on $\,S\,$, flat over $\,\La\,$. 
\item   Such a family is called a \emph{family of (maximal)  \Cm s} (respectively \emph{vector bundles})
 if for every finite dimensional
 $\,\La$-module $\,L\,$ the sheaf of $\,\kO_S$-modules $\,\kF\*_\La L\,$ is a sheaf of \Com\
 $\,\kO_S$-modules (respectively a locally free sheaf of $\,\kO$-modules).\\
 (Evidently, these properties have only to be verified for simple $\,\La$-modules.)
\item   We denote $\,\kF(L)=\kF\*_\La L\,$.
\item   We say that a family $\,\kF\,$ is \emph{generically generated over
 $\,S\,$} if there is an open subset $\,U\subseteq S\,$ such that
 $\,S\=U\,$ is a finite set of closed points and the restriction of $\,\kF\,$ onto $\,U\,$
 is generated by the image of $\,\Ga(S,\kF)\,$ in $\,\Ga(U,\kF)\,$. 
\end{itemize}
\end{defin}

 In the commutative situation we can also globalize the latter definition.

\begin{defin}
\label{fam1}
 Let $\,S,T\,$ be $\,\fK$-schemes, $\,\pr\,$ be the projection of
 $\,S\times T\,$ onto $\,T\,$ and $\,g_t\,$, where $\,t\,$ is a closed
 point of $\,T\,$, be the embedding $\,S\iso S\times t\to S\times T\,$.
\begin{itemize}
 \item
 A \emph{family of $\,\kO_S$-modules based on $\,T\,$} is, by definition, a coherent sheaf over
 $\,S\times T$ flat over $\,T\,$.
 \item
 If $\,\kM\,$ is a family of sheaves over $\,S\,$ based on $\,T\,$ and
 $\,t\,$ is a closed point of $\,T\,$, denote by $\,\kM(t)\,$ the sheaf
 $\,g_t^*\kM\,$.
 \item
  Such a family is called a \emph{family
 of vector bundles} (respectively \emph{of (maximal) \Cm s}) if
 $\,\kM(t)\,$ is locally free (respectively maximal \Cm) for each closed
 point $\,t\in T\,$.
 \item
  We say that a family $\,\kF\,$ is \emph{generically generated over
 $\,S\,$} if there is an open subset $\,U\subseteq S\,$ such that
 $\,S\=U\,$ is a finite set of closed points and the restriction onto
 $\,U\times T\,$ of the natural homomorphism $\,\pr^*\pr_*\kF\to\kF\,$ is
 an epimorphism.
\end{itemize}
\end{defin}

 General properties of flat families \cite[Sections 6.1--6.3]{gr} imply the following.

\begin{prop}
\label{fam2}
\begin{enumerate}
\item  Suppose that $\,T\,$ is a \Com{} $\,\fK$-scheme. A
 coherent sheaf $\,\kM\,$ on $\,S\times T\,$, flat over $\,T\,$, is a
 family of (maximal) \Cm s \iff $\,\kM\,$ is a sheaf of (maximal)
 \Com\ $\,\kO_{S\times T}$-modules. 
\item      Suppose that $\,T\,$ is a smooth $\,\fK$-scheme. A
 coherent sheaf $\,\kM\,$ on $\,S\times T\,$, flat over $\,T\,$, is a
 family of vector bundles \iff $\,\kM\,$ is a vector bundle over
 $\,S\times T\,$ (i.e. a locally free sheaf of $\,\kO_{S\times T}$-modules).
\end{enumerate}
\end{prop}
		
 We return to \DN and prove some results analogous to those of Kahn
 about the relations between \Cm s and vector bundles. In this case we
 call families of sheaves on $\,X=\supseteqc A\,$ \emph{families of
 $\,A$-modules}. Since we deal mostly with affine or even non-commutative base, we
 restrict our considerations to this case. The globalization of the obtained results is more
 or less evident.

 A family of vector bundles over $\,\tX\,$ based on an algebra $\,\La\,$ is called \emph{full} if it is isomorphic to
 $\,(\pi^*\kM)\dch\,$, where $\,\kM\,$ is a family of \Com{} $\,A$-modules based on $\,\La\,$. 

 Recall that an algebra 
 $\,\La\,$ is said to be \emph{hereditary} if $\,\gdim\La\le1\,$ (both left and right). Especially the local rings of points
 of any smooth curve, as well as free non-commutative algebras, and, more generally, path algebras
 of (oriented) graphs (cf. \cite[Section III.6]{dk}) are always hereditary.

\begin{prop}
\label{kahn3}\emph{(Cf.~\cite[Proposition 1.2]{kah}.)\ }
 Suppose that a family $\,\kF\,$ of vector bundles over $\,\tX\,$
 based on a hereditary algebra $\,\La\,$ (or on a smooth curve) satisfies the following conditions: 
 \begin{enumerate}
  \item
  $\,\kF\,$ is generically generated over $\,\tX\,$.
  \item
  The restriction map $\,\rH^0(\tX,\kF)\to\rH^0(\pX,\kF)\,$ is an epimorphism
 (or, equivalently, the map $\,\rH^1_E(\tX,\kF)\to\rH^1(\tX,\kF)\,$ is a monomorphism). 
 \end{enumerate}
 Then $\,\kF\,$ is full.
\end{prop}
\begin{proof}
 The claim is evidently local, so we only have to prove it in the case of algebras.
 Set $\,\kM=\pi_*\kF\,$, $\,\kF'=\pi^*\kM/(\mathrm{torsion})\,$ and consider $\,\kF'\,$
 as a subsheaf of $\,\kF\,$ in a natural way. Note that if $\,U\subseteq X\,$ is an open subset
 and $\,\pi^{-1}(U)=\bigcup U_i\,$ is an affine covering of $\,\pi^{-1}(U)\,$, $\,\kM(U)\,$ is a submodule
 of $\,\bigoplus_i\kF(U_i)\,$. Since $\,\La\,$ is hereditary, a submodule of a
 flat $\,\La$-module is flat, so $\,\kM\,$ is flat over $\,\La\,$. It is coherent since $\,\pi\,$ is proper
 and $\,\kM(L)\iso\pi_*\kF(L)\,$ is \Com\ by \cite[Proposition 6.3.1]{gr}, so $\,\kM\,$ is a family of
 \Cm s based on $\,\La\,$. Moreover, since $\,\kF\,$ is generically generated over
 $\,\tX\,$, $\,\cod\supp\kF/\kF'\ge2\,$. But $\,\tX\,$ is normal,
 hence any sheaf of the form $\,\kN^\vee\,$ is completely defined by
 its stalks at the points of height $1$. In particular,
 $\,(\pi^*\kM)^\vee\iso(\kF')^\vee\iso\kF^\vee\,$, so
 $\,\kF\iso\kF\dch\iso(\pi^*\kM)\dch\,$.
\end{proof}

\begin{prop}
\label{kahn4}\emph{(Cf. \cite[Proposition 1.6]{kah}.)\ }
 Suppose that $\,Z\,$ is a weak reduction cycle on $\,\tX\,$,
 $\,\kF\,$ is a family of vector bundles on $\,\tX\,$ based on a hereditary algebra
 $\,\La\,$ (or on a smooth curve) and $\,\kV=\kF\*_{\kO_{\tX}}\kO_Z\,$. If
 $\,\kV\,$ is generically generated by global sections and the map
 $\,\rH^0(E,\kV(Z))\to\rH^1(E,\kV)\,$ induced by the exact sequence
$$
 0\larr \kV(-Z)\larr \kV_2\larr \kV\larr 0\,,
$$
 where $\,\kV_2=\kF\*_{\kO_{\tX}}\kO_{2Z}\,$, is a monomorphism, then
 $\,\kF\,$ is full. 
\end{prop}
\begin{proof}
 Again one only has to consider the case of algebras; then one simply has to
 follow literally the proof of Kahn using proposition~\ref{kahn3}.
\end{proof}

\begin{prop}
\label{kahn5}\emph{(Cf.~\cite[Proposition 1.9]{kah}.)\ }
 Suppose that $\,Z\,$ is a weak reduction cycle on $\,\tX\,$ and a
 family $\,\kV\,$ of vector bundles over $\,Z\,$ based on a hereditary 
 algebra $\,\La\,$ (for instance, on a smooth affine curve) satisfies the following conditions: 
 \begin{enumerate}
  \item
  $\,\kV\,$ is generically generated over $\,Z\,$.
  \item
  There exists an extension of $\,\kV\,$ to a family of vector bundles
 $\,\kV_2\,$ on $\,2Z\,$ such that the induced map $\,\rH^0(E,\kV(Z))\to\rH^1(E,\kV)\,$
 is injective.
 \end{enumerate}
 Then there is a full family of vector bundles $\,\kF\,$ over
 $\,\tX\,$ such that $\,\kV\iso\kF\*_{\kO_{tX}}\kO_Z\,$. If $\,\La\,$ is finite dimensional,
 such a family is unique.
\end{prop}
\begin{proof}
 Since $\,\rH^2(E,\kN)=0\,$ for each coherent sheaf over
 $\,Z\,$, any family of vector bundles over $\,2Z\,$ can be
 lifted to a family of vector bundles over $\,nZ\,$ for every $\,n\,$,
 hence to a family of vector bundles over $\,\tX\,$. 
 Proposition~\ref{kahn4} implies that $\,\kF\,$ is full. If $\,\La\,$ is finite dimensional,
 the last assertion follows directly from \cite[Proposition 1.9]{kah}.
\end{proof}

 For minimally elliptic singularities one gets a simpler version.

\begin{prop}
\label{kahn6}
\emph{(Cf. \cite[Theorem 2.1]{kah}.)\ }
 Let $\,\pi:\tX\to X\,$ be the minimal resolution of a minimally
 elliptic singularity and $\,Z\,$ be the fundamental cycle on
 $\,\tX\,$. Suppose that a family of vector bundles $\,\kG\,$ on
 $\,Z\,$ based on a hereditary algebra $\,\La\,$ (for instance, on a smooth affine curve)
 satisfies the following conditions:
 \begin{enumerate}
  \item
  $\,\kG\,$ is generically generated over $\,Z\,$.
  \item
  $\,\rH^1(E,\kG)=0\,$.
  \item
  $\,\rH^0(E,\kG(Z))=\kP\,$ is flat (or, equivalently, projective) as $\,\La$-module.
 \end{enumerate}
 Then there is a full family of vector bundles $\,\kF\,$ over
 $\,\tX\,$ such that $\,\kF\*_{\kO_{\tX}}\kO_Z\iso\kG\+\kO_Z\*\kP\,$. If $\,\La\,$ is finite dimensional,
 such a family is unique.
\end{prop}
\begin{proof}
 In this situation we can again repeat the arguments of the corresponding proof from
 \cite[Page 151]{kah} to show that the family $\,\kV=\kG\+\kO_Z\*\kP\,$ satisfies the
 conditions of Proposition~\ref{kahn5}. 
\end{proof}

\begin{remk}\label{glob}
 Suppose that in Propositions~\ref{kahn5} and~\ref{kahn6} we consider families based on a smooth affine curve $\,T\,$.
 Then the uniqueness assertion from \cite[Proposition 1.9]{kah} can be applied to the generic fibres. It implies
 that for any two liftings $\,\kF,\,\kF'\,$ there is an open dense subset $\,U\subseteq T\,$ such that $\,\kF|_U\iso
 \kF'|_U\,$. On the other hand, if the base $\,T\,$ is a projective smooth curve, we can in general only lift an
 appropriate family $\,\kV\,$ of vector bundles over $\,Z\,$ to a family of vector bundles over $\,\tX\,$
 on an open subset $\,U\subseteq T\,$. Moreover, on a smaller open subset such a lifting is unique. This consideration implies
 that there is an etale covering $\,\th:T'\to T\,$ such that $\,\th^*\kV\,$ can be lifted to a full family of vector bundles
 over $\,\tX\,$, hence gives rise to a family of \Com\ $\,A$-modules with the base $\,T'\,$.
\end{remk}

 For cusp singularities we can now precise Corollary~\ref{cusp1} in
 the following way. (We use the notations of Section~\ref{sec2}.)
\begin{corol}
\label{cusp2}
 Let $\,X=\supseteqc A\,$ be a cusp singularity. For each aperiodic
 $\,s$-sequence $\,\fD\ge0\,$ there is a family of \Com{} $\,A$-modules
 $\,\kM(\fD)\,$ based on $\,T\,$, where $\,T=\fK^*\=\set{1}\,$ if
 $\,\fD=\oO\,$ or $\,\fD=\fB\,$, and $\,T=\fK^*\,$ otherwise, such that
 $\,M(\fD,m,\la)\iso \kM\*_{\kO_T}L(m,\la)\,$ for every $\,m\in\mN\,$
 and every $\,\la\in T\,$, where $\,L(m,\la)=\kO_T/\gM_\la^m\,$ ($\,\gM_\la\,$
 is the maximal ideal corresponding to the point $\,\la\in T\,$).
\end{corol}
\begin{proof}
 Consider the family of vector bundles over $\,\tE\,$:
 $\,\tG=\bigoplus_{i=1}^rs\kO_i\*\kO_T\,$. Define $\,\kG(\fD)\,$ as the
 preimage in $\,\tG\,$ of the sub-bundle in $\,\tG/\kJ\tG\,$ generated
 by the elements $\,\fE_i'+\fE_{i+1}''\,$ for $\,i\ne rs\,$ and the
 elements $\,\fE_{rs}'+t\fE''_1\,$, where $\,t\,$ denotes the
 coordinate on $\,T\,$ and $\,\fE_i',\,\fE_i''\,$ are defined as on
 page~\pageref{ei}, after Theorem~\ref{dg}. It is easy to see that
 $\,\kG(\fD)\*_{\kO_T}L(m,\la)\iso\kG(\fD,m,\la)\,$.
 Corollaries~\ref{cusp0} and \ref{cusp1} show that such a family
 satisfies the conditions of Proposition~\ref{kahn6}. Hence it gives a
 full family of vector bundles $\,\kF\,$ over $\,\tX\,$ such that
 $\,(\kF\*_{\kO_T}L(m,\la))\*_{\kO_{\tX}}\kO_E\iso
 mn(\fD)\kO_E\+\kG(\fD,m\la)\,$, where $\,n(\fD)=
 \sum_{i=1}^{rs}(d_i-b_i+1)-\th(\fD-\fB^r)\,$. Therefore, we obtain a
 family $\,\kM(\fD)\,$ of \Com{} $\,A$-modules just as we need. 
\end{proof}

 Note that among $\,\kM(\fD)\*_{\kO_T}L(m,\la)\,$ one obtains all
 indecomposable \Com{} $\,A$-modules except $\,A\,$ and
 $\,M(\fB,m,1)\,$. 

\section{\Com{} tame and wild singularities}
\label{sec4}

 We provide now some formal definitions of ``tame'' and ``wild''
 singularities (with respect to the classification of \Cm s),
 analogous to the usual ones from the representation theory of
 algebras. Again we consider the case of $\,\fK$-algebras. 

\begin{defin}\label{cmtame}
\begin{enumerate}
\item 
 We call a surface singularity $\,X=\supseteqc A\,$ \emph{\Com\ tame} if it is not
 \Com\ finite and there is a set $\,\dM=\set{\kM^\al}\,$ of families of \Com\
 $\,A$-modules such that
 \begin{itemize}
  \item  Each $\,\kM^\al\,$ is based on a smooth connected curve $\,C^\al\,$.
  \item  The set $\,\dM_r=\setsuch{\kM^\al}{\rk\kM^\al=r}\,$ is finite; we denote the number of 
 elements in it by $\,d(\dM,r)\,$.
  \item  Almost all indecomposable \Com{} $\,A$-modules of rank $\,r\,$ (i.e. all except 
 finitely many isomorphism classes) are isomorphic to $\,\kM^\al(c)\,$ for some $\,\kM^\al\in\dM_r\,$
 and some closed point $\,c\in C^\al\,$.
 \end{itemize}
 We call a set $\,\dM\,$ with these properties a \emph{parametrizing set} for \Com\ $\, A$-modules.
 We also denote by $\,d(r,X)=\min\set{d(\dM,r)}\,$, where $\,\dM\,$ runs through all parametrizing
 sets for \Com\ $\,A$-modules.
 \item  We call a tame surface singularity \emph{bounded} (or of \emph{polynomial growth}) if there
 is a polynomial $\,\vi(r)\,$ such that $\,d(r,X)\le\vi(r)\,$ for all $\,r\,$, and \emph{unbounded} otherwise.
 We say that this singularity is \emph{of exponential growth} if $\,d(r,X)\,$ growths exponentially when
 $\,r\to\8\,$.
 \end{enumerate}
 \end{defin}

\begin{exam}
\label{extame}
\begin{enumerate}
 \item
  Corollary~\ref{cusp2} shows immediately that every cusp singularity
  is \Com{} tame. Moreover, it follows from \cite{dg2} that it is of exponential growth.
 \item
  Just in the same way, the results of \cite{kah} show that every
  simple elliptic singularity is also \Com{} tame, bounded (the corresponding
  families are based either on $\,\pic^0E\,$ or on
  $\,\pic^0E\=\set{0}\,$). 
\end{enumerate} 
\end{exam}

\begin{remk}
 Denote by $\,\pic^0E\,$ the subgroup of the Picard group of a
 (singular) curve $\,E\,$ consisting of the classes of all line
 bundles $\,\kL\,$ such that the degree of the restriction of
 $\,\kL\,$ onto every irreducible component of $\,E\,$ equals $0$. One
 can easily see that $\,\pic^0E\iso\fK^*\,$ for  every cusp
 singularity. Since the parametrizing families for cusp singularities are based either on
 $\,\fK^*\,$ or on $\,\fK^*\=\set1\,$, both examples above have a lot of common features. 
\end{remk}

\begin{defin}
\label{cmwild}
\begin{enumerate}
\item   A family $\,\kM\,$ of  \Com\ $\, A$-modules based on an algebra
 $\,\La\,$ is called \emph{strict} if for every finite dimensional $\,\La$-modules $\,L,L'\,$
  \begin{itemize}
  \item  $\,\kM(L)\iso\kM(L')\,$ \iff $\,L\iso L'\,$;
  \item   $\,\kM(L)\,$ is indecomposable \iff so is $\,L\,$.
  \end{itemize}
 \item  We call a surface singularity $\,X=\supseteqc A\,$ \emph{\Com{} wild} if for every
 finitely generated (not necessarily commutative) $\,\fK$-algebra $\,\La\,$ there is a
 strict family of \Com\ $\,A$-modules based on $\,\La\,$.
\end{enumerate}
\end{defin}

\begin{remk}
\label{re1}
 It is well known \cite{dr1} that to prove that $\,X\,$ is \Com{} wild one
 only has to find a strict family for one of the following algebras $\,\La\,$:
 \begin{itemize}
  \item[(i)]
  $\,\fK\gnr{x,y}\,$, the free non-commutative algebra with $2$
 generators;
  \item[(ii)]
  $\,\fK[\,x,y\,]\,$, the polynomial algebra with $2$ generators;
  \item[(iii)]\label{ca3}
  $\,\fK[[\,x,y\,]]\,$, the power series algebra with $2$ generators. 
 \item[(iv)]\label{ca4}
  $\,\fK\Ga_3\,$, the path algebra of the graph
$$
  \Ga_3 = \xymatrix{ \bullet \ar@<.75ex>[rr] \ar[rr] \ar@<-.75ex>[rr]&&\bullet }								
$$
 It is the 5-dimensional algebra with a basis $\,\set{e_1,e_2,a_1,a_2,a_3}\,$ and the multiplication:
 $\,e_i^2=e_i,\,a_ie_1=e_2a_i=a_i\,$, all other products are zero (cf. for instance \cite{dk}).
 \end{itemize}
 Cases (i) and (iv) are the most appropriate for our purpose since they are hereditary
 (and (iv) is even finite dimensional), so we can use all results on families from the preceding section.
 Since we do not know whether those results do also hold for non-hereditary algebras, for instance, for
 $\,\fK[x,y]\,$ (although it may be conjectured), the proof of our main results rests indeed upon using families
 with non-commutative bases.

 Obviously, if $\,X\,$ is \Com{} wild it has families $\,\kM\,$ of \Cm s
 based on any given algebraic variety $\,T\,$, such that all modules 
 $\,\kM(t)\,$ with $\,t\in T\,$ are indecomposable and pairwise non-isomorphic.
\end{remk}

\begin{theorem}
\label{minel}
 Suppose that a minimally elliptic singularity  $\,X\,$ is neither simple
 elliptic nor cusp. Then it is \Com{} wild.
\end{theorem}
\begin{proof}
 Recall \cite{dg2} that $\,E\,$ is vector bundle wild if $\,E\,$ is
 neither a smooth rational or elliptic curve, nor a cyclic or linear
 configuration, where a \emph{linear configuration} is a curve $\,E\,$
 such that all its components $\,E_i\,$ ($\,i=1,\dots,s\,$) are smooth
 rational, $\,E_i\cap E_j=0\,$ if $\,j\ne i\pm1\,$, while $\,E_i\,$
 and $\,E_{i+1}\,$ intersect transversally in exactly one point for
 $\,i=1,\dots,s-1\,$ (in other words, its dual graph is of type $\rA_n$).
 The latter case is impossible since $\,X\,$ is not
 rational. So, if $\,X\,$ is neither simple elliptic nor a cusp, the
 fundamental cycle $\,Z\,$ is vector bundle wild, i.e. there is a strict
 family of vector bundles $\,\kG\,$ over $\,Z\,$ based on
  $\,\fK\Ga_3\,$. By Serre's Theorems \cite[Theorems II.5.17 and
 III.5.2]{har}, there is an integer $\,n\,$ such that $\,\kG(n)\,$ is
 generated by global sections and $\,\rH^1(E,\kG(n))=0\,$. Since
 $\,\kG(n)\,$ is obviously also strict, one can suppose that $\,\kG\,$ 
 itself has these properties.  By Proposition~\ref{kahn6}, there is a full
 family $\,\kF\,$ of vector bundles over $\,\tX\,$ based on $\,\fK\Ga_3\,$
 such that $\,\kF\*_{\kO_{\tX}}\kO_Z\iso\kG\+\kO_Z\*\pr_*\kG\,$. Let
 $\,\kF=(\pi^*\kM)\dch\,$, where $\,\kM\,$ is a family of \Cm s. Since
 $\,\kG\,$ is strict, Kahn's theorem (Theorem~\ref{kahn2}) implies that
 $\,\kM\,$ is also strict, hence $\,X\,$ is \Com{} wild.
\end{proof}

 Consider the following procedure that allows to obtain new \Com{} tame singularities.

 Let $\,B\,$ be a local
 noetherian ring with maximal ideal $\,\gM\,$, $\,A\supseteq B\,$ be a finite $\,B$-algebra. We call this extension
 (or the corresponding morphism $\,\supseteqc A\to \supseteqc B\,$)
 \emph{split} if the embedding $\,B\to A\,$ splits as monomorphism of $\,B$-modules, i.e. $\,A\iso B\+(A/B)\,$.
 It is called \emph{unramified} if $\,A/\gM A\,$ is a separable $\,B/\gM$-algebra, or, equivalently,
 the natural epimorphism $\,\eps:A\*_BA\to A\,$ of $\,A$-bimodules splits \cite{ab}.
 $\,A\,$ is called \emph{unramified in codimension 1} if the extension $\,B_\gP\subseteq A_\gP\,$ is
 unramified for every prime ideal $\,\gP\subseteq B\,$ of height 1. We also denote by $\,A\bt M\,$, where $\,M\,$
 is a $\,B$-module, the second dual $\,(A\*_BM)\dch\,$ of the $\,A$-module $\,A\*_BM\,$.

  \begin{lemma}\label{exten}
 Let $\,B\subseteq A\,$ be a finite extension of normal rings.
 This extension is unramified in codimension 1 \iff the epimorphism of $\,A$-bimodules
 $\,\eps\dch:A\bt A\to A\,$ splits. 
\end{lemma} 
\begin{proof}
 If $\,\eps\dch\,$ splits, $\,\eps_\gP\dch\,$ also splits for every prime $\,\gP\,$. If $\,\he\gP=1\,$,
 $\,\eps\dch_\gP\,$ coincides with $\,\eps_\gP\,$, thus the extension $\,B_\gP\subseteq A_\gP\,$ is unramified.
 Suppose now that $\,B\subseteq A\,$ is unramified in codimension 1.  Then the extension
 $\,K\subseteq L\,$ is separable, hence the exact sequence 
$$
   0\larr J\larr L\*_KL\stackrel\phi\larr L\larr0  
$$
 splits and $\,L\*_K L\,$ is a semi-simple ring, so $\,L\*_KL=J\+J'\,$, where $\,J'=\mathrm{Ann}_{L\*_KL}J\,$,
 $\,\phi(J')=L\,$, and the restriction of $\,\phi\,$ onto $\,J'\,$ is an isomorphism.
 Obviously, $\,A\bt A\,$ is a $\,B$-subalgebra in $\,L\*_KL\,$ and $\,\eps\dch\,$ is the restriction of $\,\phi\,$
 onto $\,A\bt A\,$. Set $\,I=J\cap(A\bt A),\, I'=J'\cap(A\bt A)\,$. Then $\,\phi\,$ induces an epimorphism of $\,L$-bimodules
 $\,\phi_*:\Hom_{L\mbox{-}L}(L,L\*_KL)\to\Hom_{L\mbox{-}L}(L,L)\,$. Since $\,L\iso (L\*_KL)/J\,$, $\,\phi_*\,$
 can be identified with a mapping $\,J'\to L\,$, so it is an isomorphism. Moreover, $\,\phi_*(I'_\gP)=A_\gP\,$ for every
 prime ideal $\,\gP\,$ of height 1, as the extension $\,B_\gP\subset A_\gP\,$ is unramified.
 But obviously $\,I'=\bigcap_{\he\gP=1}I'_\gP\,$ and $\,A=\cap_{\he\gP=1}A_\gP\,$, so $\,\phi_*(I')=A\,$ and
 $\,A\bt A\to A\,$ splits.
\end{proof}

\begin{prop}
\label{ext1}
 Let $\,X=\supseteqc A,\ Y=\supseteqc B\,$ be normal surface singularities and
 $\,X\to Y\,$ be a finite surjective morphism given by an extension $\,B\to A\,$.
 \begin{enumerate}
  \item
  If $\,X\to Y\,$ is split and $\,X\,$ is \Com{} tame, so is also $\,Y$; moreover, if $\,X\,$ is
 bounded, so is $\,Y\,$.
  \item
  If $\,X\to Y\,$ is unramified in codimension 1 and $\,Y\,$ is \Com{} tame, so is also $\,X$;
 moreover, if $\,Y\,$ is bounded, so is $\,X\,$.
 \end{enumerate}
\end{prop}
\begin{proof}
 (1) Suppose that this extension is split and $\,\rk_BA=m\,$. If $\,N\,$ is a \Com{}
    $\,B$-module of rank $\,r\,$, it is a direct summand of $\,A\*_BN\,$, hence
    of $\,A\bt N\,$ as $\,B$-module. Note that $\,\rk_AA\bt N=r\,$ too, so $\,\rk_BA\bt N=mr\,$.
   Suppose that $\,A\,$ is \Com{} tame. For every family $\,\kM^\al\,$ from
    Definition~\ref{cmtame} let $\,t^\al\,$ be the general point of
    the curve $\,C^\al\,$ and $\,\kM^\al(t^\al)=\bigoplus_\be\tN^{\al\be}\,$, where
   $\,\tN^{\al\be}\,$ are indecomposable $\,B\*\fK(t^\al)$-modules. There is an open
    subset $\,U^\al\subseteq C^\al\,$ and a decomposition
    $\,\kM^\al|_{U^\al}\iso\bigoplus_\be\kN^{\al\be}\,$ (as family of
    $\,B$-modules) such that $\,\kN^{\al\be}(c)\,$ is indecomposable
    for all $\,c\in U^\al\,$. Denote by $\,\dN_r\,$ the set of all $\,\kN^{\al\be}\,$ such that 
   $\,\rk\kN^{\al\be}=r\,$, while $\,r/m\le\rk\kM^\al\le r\,$. It is a finite set having at most
  $\,\sum_{j=[r/m]}^r[mj/r]d(\dM,j)\,$ elements. We also denote by $\,\dM^*_s\,$ the set of indecomposable
  \Com\ $\,A$-modules of rank $\,s\,$ that are not of the form $\,\kM^\al(c)\ (c\in U_\al)\,$ and by $\,\dN^*_r\,$
 the set of indecomposable $\,B$-modules of rank $\,r\,$ that are direct summands (over $\,B$) of modules
 $\,M\in\dM^*_s\,$ with $\,r/m\le s\le r\,$. Certainly $\,\dM^*_s\,$ and $\,\dN^*_r\,$ are both finite.

    Let $\,N\,$ be any indecomposable \Com{} $\,B$-module of rank $\,r\,$. It is a
    direct summand of an indecomposable \Com{} $\,A$-module
    $\,M\,$ of rank $\,s\,$ with $\,r/m\le s\le r\,$. If $\,\rk M=s\,$, either $\,M\in\dM^*_s\,$ or
 $\,M\iso\kM^\al(c)\,$ for some $\,\kM_\al\in\dM_s\,$ and some $\,c\in C^\al\,$. In the former case
 $\,N\in\dN^*_r\,$, while in the latter case $\,N\iso\kN^{\al\be}(c)\,$ for some $\,\kN^{\al\be}\in\dN_r\,$
 and some $\,c\in U^\al\,$. Therefore $\,Y\,$ is \Com\ tame. Moreover, $\,d(r,Y)\le\sum_{j=[r/m]}^r[mj/r]d(j,X)\,$,
 hence if $\,X\,$ is bounded, so is $\,Y\,$.

\smallskip
 (2) Suppose that this extension is unramified in codimension 1. If $\,M\,$ is any
    \Com{} $\,A$-module, one easily verifies that $\,(A\bt M)^\vee\iso
    ((A\bt A)\*_AM)^\vee\,$. Hence $\,M^\vee\,$ is isomorphic to a
    direct summand of $\,(A\bt M)^\vee\,$ and $\,M\,$ is isomorphic to
    a direct summand of $\,A\bt M\,$. So, there is an indecomposable \Com{} $\,B$-module
  $\,N\,$ such that $\,M\,$ is isomorphic to a direct summand of $\,A\bt N\,$ and
  $\,\rk_BN\le\rk_AM\le m\rk_BN\,$. Now the same considerations as above
 show that if $\,B\,$ is \Com{} tame (bounded), so is $\,A\,$.
\end{proof}

\begin{remk}\label{expon}
 Suppose that $\,X\to Y\,$ is both split and unramified in codimension 1. Then the proof of
 proposition~\ref{ext1} implies that if $\,X\,$ is of exponential growth, so is $\,Y\,$. 
\end{remk}

\smallskip
 Important examples arise from group actions.
 
 \begin{prop}\label{group}
 Let $\,G\,$ be a finite group of automorphisms of a normal surface singularity $\,A\,$,
 and let $\,B=A^G\,$ be the subalgebra of $\,G$-invariants. For a prime ideal $\,\gP\subset A\,$ set
$$
   G_\gP=\setsuch{g\in G}{g\gP=\gP\,\text{\emph{ and }}\,g\,\text{\emph{ acts trivially on }}
  \,A/\gP}.
$$
 \begin{enumerate}
\item  If $\,\chr\fK\,$ does not divide $\,n=\card G\,$, the extension $\,B\subseteq A\,$ splits.
\item  If $\,G_\gP=\set{1}\,$ for every prime ideal $\,\gP\subset A\,$ of height 1, this extension is unramified in
 codimension 1. \\
 \em In the latter case we say that $\,G\,$ \emph{acts freely in codimension} 1.
\end{enumerate}
\end{prop}
\begin{proof}
 (1) The mapping $\,A\to B,$ $\,a\mapsto\frac1n\sum_{g\in G}ga\,$, splits this extension.

 (2) Let $\,\gQ\subset B\,$ be any prime ideal of height 1, $\,\lst\gP s\,$ be all minimal prime ideals of $\,A\,$ containing
 $\,\gQ\,$. All of them are of height 1, and $\,G\,$ acts transitively on the set $\,\set{\lst\gP s}\,$ \cite[Section III.3]{ser}.
 Denote $\,K=B_\gQ/\gQ B_\gQ\,$, $\,L=A_\gQ/\gQ A_\gQ\,$ and $\,L_i=A_{\gP_i}/\gP_iA_{\gP_i}\,$.
 Since $\,\rk_BA=\card G=n\,$, $\,L\,$ is an $\,n$-dimensional $\,K$-algebra, and $\,\prod_{i=1}^sL_i\,$ is a
 factor-algebra of $\,L\,$. Set $\,G_i=\setsuch{g\in G}{g\gP=\gP}\,$. These subgroups are all of the same cardinality
 $\,m\,$ and $\,ms=n\,$. The subgroup $\,G_i\,$ acts as a group of automorphisms of $\,L_i\,$ over $\,K\,$, so
 $\,m\le\dim_KL_i\,$ and $\,n=ms\le\dim_K\prod_{i=1}^sL_i\,$. Hence, $\,L=\prod_{i=1}^sL_i\,$ and $\,m=\dim_KL_i\,$,
 so $\,L_i\,$ is a Galois extension of $\,K\,$ and $\,L\,$ is separable over $\,K\,$.
\end{proof}

\begin{corol}\label{gr1}
 \emph{We keep the notations of the preceding proposition.}
\begin{enumerate}
\item  If $\,\chr\fK\,$ does not divide $\,\card G\,$ and $\,A\,$ is \Com\ tame (bounded),
 so is $\,B\,$.
\item  If $\,G\,$ acts freely in codimension 1 and $\,B\,$ is tame (bounded), so is $\,A\,$.
\item  If both conditions hold, $\,A\,$ is of exponential growth \iff so is $\,B\,$.
\end{enumerate} 
\end{corol}

 Call a surface singularity $\,Y\,$ a \emph{simple elliptic-quotient}
 (respectively a \emph{cusp-quotient}) if it is a quotient of a simple
 elliptic singularity (respectively of a cusp) by a finite group $\,G\,$ such
 that $\,\chr\fK\,$ does not divide $\,\card G\,$. Both simple
 elliptic-quotient and cusp-quotient singularities will be called
 \emph{elliptic-quotient} singularities. If $\,\chr\fK=0\,$,
 elliptic-quotient singularities coincide with those
 \emph{log-canonical} ones which are not quotient singularities
 \cite{kaw}.

\begin{corol}
\label{log}
 Every elliptic-quotient surface singularity is \Com{} tame. Among them, simple elliptic-quotient are
 bounded, while cusp-quotient are of exponential growth.
\end{corol}

 A special case is that of \emph{Q-Gorenstein} singularities and their Gorenstein coverings
 defined as follows.
\begin{defin}\label{QG}
 Let $\,\om=\om_B\,$ be a dualizing ideal of a normal singularity $\,B\,$. Denote by $\,\om\hh k=
 \bigcap_{\he\gP=1}\om^k_\gP \ (k\in\mZ)\,$ (note that each $\,\om_\gP\,$ is a principal ideal, since $\,B_\gP\,$
 is a \dvr).
 Call $\,B\,$ \emph{Q-Gorenstein} if there is $\,n>0\,$, prime to $\,\chr\fK\,$, such that the ideal $\,\om\hh n\,$
 is principal.
\end{defin}

\begin{prop}\label{QGcov}
 Suppose that $\,B\,$ is Q-Gorenstein. Let $\,n\,$ be the smallest positive integer such that $\,\om\hh n\,$ is principal,
 $\,\om\hh n=\th B\,$. Denote $\,A=\bigoplus_{k=0}^{n-1}\om\hh k\,$ and consider it as $\,B$-algebra by setting $\,a\cdot b=
 ab/\th\,$ for $\,k+l\ge n,\,a\in\om\hh k,\,b\in\om\hh l\,$. Then $\,A\,$ is a normal Gorenstein singularity and
 $\,B=A^G\,$, where $\,G\,$ is a cyclic group of order $\,n\,$ that acts on $\,A\,$ freely in codimension 1.
 \em We call $\,A\,$ the \emph{Gorenstein covering} of $\,B\,$.
\end{prop}
\begin{proof}
  If $\,\gP\subset B\,$ is a prime ideal of height 1, $\,\om_\gP=\ga B_\gP\,$ for some $\,\ga\,$ and $\,\th=\ga^n\ze\,$
 for an invertible element $\,\ze\in B_\gP\,$. Then $\,A_\gP\iso B_\gP[t]/(t^n-\ze)\,$, so it is unramified in codimension
 1, especially $\,A_\gP\,$ is normal. Since $\,A\,$ is \Com, $\,A=\bigcap_{\he\gP=1}A_\gP\,$, so $\,A\,$ is normal itself.
  Moreover, $\,\Hom_B(\om\hh k,\om)\iso\om\hh{1-k}\,$, so $\,\om_A\iso\Hom_B(A,\om_B)\iso A\,$, thus $\,A\,$ is
 Gorenstein. If $\,K\,$ is the field of fractions of $\,B\,$, $\,L\,$ is the field of fractions of $\,A\,$, then $\,L\iso K[t]/(t^n-\ze)\,$
 is a Galois extension of $\,K\,$ with cyclic Galois group of order $\,n\,$. Therefore $\,B=A^G\,$ and $\,G\,$ acts
 freely in codimension 1.
\end{proof}

 \begin{prop}\label{QGtame}
 Let $\,B\,$ be a Q-Gorenstein surface singularity, $\,A\,$ be its Gorenstein covering. If $\,A\,$ is \Com\ tame
 (bounded, of exponential growth), so is $\,B\,$ and vice versa.
\end{prop}

\begin{erem}
 If $\,\chr\fK=0\,$, the log-canonical singularities are just Q-Gorenstein, such that their Gorenstein coverings are either
 rational double points (in the case of quotient singularities), or simple elliptic, or cusp singularities.
\end{erem}

 We call a normal surface singularity \emph{Q-elliptic} if its is Q-Gorenstein and its Gorenstein covering is
 minimally elliptic.
 \begin{corol}\label{Qel}
 A Q-elliptic singularity is \Com\ tame \iff it is elliptic-quotient; otherwise it is \Com\ wild.
\end{corol}

\section{Curve singularities $\,\rm{T_{pq}}\,$}
\label{sec5}

 Important examples of cusp singularities for $\,\chr\fK=0\,$ are the
 ``serial'' unimodal singularities \cite{arn} $\,\rT_{pqr}\,$:
 $\,A=\fK[[x,y,z]]/(x^p+y^q+z^r+\al xyz)\,$, where $\,r\le p\le q\,$
 and $\,1/p+1/q+1/r<1,\ \al\in\fK^*\,$ (in this case all values of $\,\al\,$ lead to
 isomorphic algebras). Note that if $\,1/p+1/q+1/r=1\,$, the
 corresponding singularity is simple elliptic except for finitely many special
 values of $\,\al\,$. If $\,r=2\,$, Cohen--Macaulay modules
 over this singularity are closely related to those over the
 \emph{curve singularity} $\,\rT_{pq}\,$. Indeed, in this case the
 singularity $\,\rT_{pq2}\,$ can be rewritten in the form:
 $\,A=\fK[[x,y,z]]/(z^2+x^p+y^q+\be x^2y^2)\,$ for
 $\,\be=-\al^2/4\,$. Therefore, one can use the \emph{Kn\"orrer's
 correspondence} \cite{kno,yos} described in the following proposition
 to relate \Cm s over $\,A\,$ and over the curve singularity
 $\,A/(z^2)=\fK[[x,y]]/(x^p+y^q+\be x^2y^2)\,$, denoted $\,\rT_{pq}\,$.
\begin{prop}
\label{knor}
 Let $\,f\in\fK[[\lst xn]]\,$ be a non-invertible element,
 $\,A=\fK[[\lsto xn]]/(x_0^2+f)\,$ and $\,A'=\fK[[\lst
 xn]]/(f)\,$. For every \Cm{} $\,M\,$ over $\,A\,$ let $\,\rest M=M/x_0M\,$
 and $\,M^\si\,$ be the $\,A$-module that coincide with $\,M\,$ as a
 group, while the action of the elements of $\,A\,$ is given by the
 rule: $\,a\cdot_{M^\si}v=a^\si\cdot_Mv\,$, where $\,g(\lsto
 xn)^\si=g(-x_0,\lst xn)\,$. Then: 
 \begin{itemize}
  \item
  every indecomposable \Com{} $\,A'$-module is a direct summand of
 $\,\rest M\,$ 
 for an indecomposable \Com{} $\,A$-module $\,M\,$ and $\,M\iso A\,$
 \iff $\,\rest M\iso A'\,$;
  \item
  if $\,M\,$ is an indecomposable \Com{} $\,A$-module such that
 $\,M\not\iso M^\si\,$, $\,\rest M\,$ is also indecomposable;
  \item
  if $\,M\not\iso A\,$ is an indecomposable \Com{} $\,A$-module such
 that $\,M\iso
 M^\si\,$, then $\,\rest M\,$ is a direct sum of two non-isomorphic
 indecomposable \Com{} $\,A'$-modules, which we denote by $\,\rest_1M\,$
 and $\,\rest_2M\,$.
  \item
  if $\,M,N\,$ are non-isomorphic indecomposable \Com{} $\,A$-modules
 and $\,M\not\iso N^\si\,$, then the indecomposable $\,A'$-modules
 obtained from $\,M\,$ and from $\,N\,$ as described above are also
 non-isomorphic.   
 \end{itemize}
\rm Obviously, always $\,\rest M\iso\rest(M^\si)\,$ and $\,A^\si\iso A\,$.
\end{prop}

 So, to describe the \Cm s over the curve singularity $\,\rT_{pq}\,$
 we have to know when $\,M^\si\iso M\,$ for a \Cm{} $\,M\,$ over the
 surface singularity $\,\rT_{pqr}\,$. The automorphism $\,g\to g^\si\,$
 induces an automorphism of the minimal resolution $\,\tX\,$, hence an
 automorphism of the exceptional curve $\,E\,$ and of the category of
 vector bundles over $\,E\,$. We denote all these automorphisms by
 $\,\si\,$ too. The following result is immediate (cf. the definition of $\,R_E\,$
 from Theorem~\ref{kahn1}).
\begin{prop}
 $\,R_EM^\si\iso(R_EM)^\si\,$ for every \Com{} $\,A$-module $\,M\,$.
\end{prop}

 From the description of the minimal resolutions of $\,\rT_{pqr}\,$
 given, for instance, in \cite{kar,lau}, one can deduce the following
 shape of the exceptional curve $\,E\,$ for the minimal resolution of
 $\,\rT_{pq2}\,$ and for the action of the automorphism $\,\si\,$ on
 it. 
\begin{prop}
\label{pq0}
{\em(We use the notations from Section~\ref{sec2}.)} 
\begin{enumerate}
 \item
  If $\,p=3\,$ (hence $\,q\ge7\,$), then $\,s=q-6\,$. If $\,q>7\,$,
  $\,E\,$ has $1$  component $\,E_1\,$ with self-intersection $\,-3\,$
  and $\,s-1\,$ components $\,E_2,\dots,E_s\,$ with self-intersection
  $\,-2\,$; $\,E_1^\si=E_1\,$ and $\,E_i^\si=E_{s+2-i}\,$. If
  $\,q=7\,$, $\,E\,$ is irreducible with one node and $\,E\cd
  E=-1\,$. In both cases $\,\si(p_1'')=p_1'\,$. 
 \item
  If $\,p=4\,$ (hence $\,q\ge5\,$), then $\,s=q-4\,$. If $\,q>5\,$,
  $\,E\,$ has $1$ component $\,E_1\,$ with  $\,E_1\cd E_1=-4\,$ and
  $\,s-1\,$ components $\,E_2,\dots,E_s\,$ with $\,E_i\cd
  E_i=-2\,$. If $\,q=5\,$, $\,E\,$ is irreducible with one node and
  $\,E\cd E=-2\,$. The action of $\,\si\,$ is the same as in the
  previous case.
 \item
  If $\,p\ge5\,$, then $\,s=p+q-8\,$. Put $\,t=p-3\,$. $\,E\,$ has
  $\,2\,$ components $\,E_1\,$ and $\,E_t\,$ with $\,E_i\cd
  E_i=-3\,$ and $\,s-2\,$ components with $\,E_i\cd E_i=-2\,$ (no one
  if $\,p=q=5\,$). $\,E_i^\si=E_{t+1-i}\,$ if
  $\,1\le s\le t\,$ and $\,E_i^\si=E_{s+t+1-i}\,$ if $\,t<i\le s\,$;
  $\,\si(p''_1)=p_t'\,$.
\end{enumerate}
\end{prop}
 Given an $\,s$-sequence $\,\fD=\row d{rs}\,$, denote by $\,\fD^\si\,$
 the sequence $\,\row{d'}{rs}\,$ obtained from $\,\fD\,$ by the
 following procedure (we keep the notations of Proposition~\ref{pq0}): 
\begin{itemize}
 \item
  if $\,p=3\,$ or $\,p=4\,$, then $\,d'_1=d_1\,$ and
  $\,d'_i=d_{rs+2-i}\,$ for $\,1<i\le rs\,$;
 \item
  if $\,p\ge5\,$, then $\,d'_i=d_{t+1-i}\,$ for $\,1\le i\le t\,$ and
  $\,d'_i=d_{rs+t+1-i}\,$ for $\,t<i\le rs\,$. 
\end{itemize}
 We call an $\,s$-sequence $\,\fD\,$ \emph{$\,\si$-symmetric} if
 $\,\fD^\si=\fD^k\,$, where $\,\fD^k\,$ is an $\,s$-shift of $\,\fD\,$
 (cf. page~\pageref{dg}). Then the description of \Cm s from
 Section~\ref{sec2} implies the following results. 
\begin{corol}
\label{pq1}
 Let $\,M=M(\fD,m,\la)\,$ be an indecomposable \Cm{} over the surface 
 singularity $\,\rT_{pq2}\,$ (cf. Corollary~\ref{cusp1}). Then
 $\,M^\si\iso M(\fD^\si,m,1/\la)\,$. In particular, $\,M\iso M^\si\,$
 \iff the sequence $\,\fD\,$ is $\,\si$-symmetric and $\,\la=\pm1\,$.
\end{corol}

\begin{corol}
\label{pq2}
 The indecomposable \Cm s over the curve singularity
 $\,A'=\fK[[x,y]]/(x^p+y^q+\be x^2y^2)\,$ of type $\,\rT_{pq}\,$
 ($\,q\ge p,\,1/p+1/q<1/2\,$) are the following:
 \begin{itemize}
  \item
  $\,N(\fD,m,\la)=\rest M(\fD,m,\la)\,$, where either the sequence
 $\,\fD\,$ is not $\,\si$-symmetric or $\,\la\ne\pm1\,$;
  \item
  $\,N_i(\fD,m,\pm1)=\rest_iM(\fD,m,\pm1)\,$ ($\,i=1,2\,$), where 
  $\,\fD\,$ is $\,\si$-symmetric and either $\,\fD\ne\oO\,$
  or $\,\la=-1\,$;
 \item
 the regular module $\,A'\,$.
 \end{itemize}
 The only isomorphisms between so defined \Cm s are
 $\,N(\fD,m,\la)\iso N(\fD^\si,m,1/\la)\,$. 
\end{corol}
 Note that the modules $\,M(\fD,m,\la)\,$ with $\,\la\ne\pm1\,$ form
 rational families based on $\,(\fK^*\=\set{1,-1})/\mZ_2\,$, where the
 cyclic group $\,\mZ_2\,$ acts on $\,\fK\,$ mapping $\,\la\,$ to
 $\,1/\la\,$. One easily sees that this factor is isomorphic to
 $\,\fK^*\=\set{1,-1}\,$ itself.

\section{Hypersurface singularities}

 We also can use the previous results and Kn\"orrer's correspondence to
 describe \Cm s on hypersurface singularities of type $\,\rT_{pqr}\,$ in any dimension.
 Such a hypersurface is given by the equation
 $\,x^p+y^q+z^r+\al xyz+\sum_{j=1}^kv_i^2\,$ for some $\,k\,$ (maybe,
 $0$). Recall the main results concerning Kn\"orrer's correspondence
 \cite{kno,yos}. 

\begin{prop}
\label{knor2}
 Let $\,f\in\fK[[\lst xn]]\,$ be a non-invertible element,
 $\,A=\fK[[\lst xn]]/(f)\,$ and $\,\tA=\fK[[\lsto
 xn]]/(f+x^2_0)\,$. For every \Cm{} $\,M\,$ over $\,A\,$ let $\,\syz
 M\,$ denote the first syzygy module of $\,M\,$ considered as
 $\,\tA$-module via the natural epimorphism $\,\tA\to
 A\iso\tA/(x_0)\,$. If $\,M\,$ has no free direct summands, let $\,\Om
 M\,$ denote its first syzygy as $\,A$-module. Then: 
 \begin{itemize}
  \item
  every indecomposable \Com{} $\,\tA$-module is a direct summand of
 $\,\syz M\,$ for an indecomposable \Com{} $\,A$-module $\,M\,$ and
 $\,M\iso A\,$ \iff $\,\syz M\iso \tA\,$; 
  \item
  if $\,M\,$ is an indecomposable \Com{} $\,A$-module such that
 $\,M\not\iso \Om M\,$, $\,\syz M\,$ is also indecomposable;
  \item
  if $\,M\not\iso A\,$ is an indecomposable \Com{} $\,A$-module such
 that $\,M\iso
 \Om M\,$, then $\,\syz M\,$ is a direct sum of two non-isomorphic
 indecomposable \Com{} $\,\tA$-modules, which we denote by $\,\syz_1M\,$
 and $\,\syz_2M\,$.
  \item
  if $\,M,N\,$ are non-isomorphic indecomposable \Com{} $\,A$-modules
 and $\,M\not\iso\Om N \,$, then the indecomposable $\,\tA$-modules
 obtained from $\,M\,$ and from $\,N\,$ as described above are also
 non-isomorphic.
  \item
  if $\,M\,$ has no free direct summands, then $\,\syz(\Om M)\iso\syz
 M\,$. 
 \end{itemize}
\end{prop}
 
 Obviously, this proposition implies immediately that if $\,A\,$ is
 \Com{} tame, so is $\,\tA\,$ (and vice versa in view of
 Proposition~\ref{knor}), just as in Proposition~\ref{ext1}. In
 particular, as we have already mentioned, for $\,k=0\,$, the (surface)
 singularity of type $\,\rT_{pqr}\,$ is simple elliptic for $\,1/p+1/q+1/r=1\,$
 and a cusp for $\,1/p+1/q+1/r<1\,$ \cite{lau}. Hence the results of Section~3 imply

\begin{corol}
 Every hypersurface singularity of type $\,\rT_{pqr}\,$ is \Com{}
 tame. 
\end{corol}

\section{Auslander--Reiten quivers}

 In this section we calculate Auslander--Reiten quivers for \Cm s and
 related vector bundles. Recall that an \AR{} in a category of modules
 or of vector bundles is a non-split exact sequence
$$
 0\larr M'\stackrel g\larr N \stackrel f\larr M\larr 0\,,
$$
 where $\,M,M'\,$ are indecomposable, such that any homomorphism
 $\,N'\to M\,$ that is not a split epimorphism factors through
 $\,f\,$ (equivalently, any homomorphism $\,M'\to N'\,$ that is not a
 split monomorphism factors through $\,g\,$). Each of the modules
 $\,M,M'\,$ uniquely defines the second one. They write $\,M'=\tau
 M\,$ and call $\,M'\,$ the \emph{\AT} of $\,M\,$. In \cite{ar1} it was
 proved that for any indecomposable vector bundle $\,\kG\,$ on a
 projective curve $\,E\,$ there is an \AR{}
$$
 0\larr \om_E\*_{\kO_E}\kG\larr \kE \larr \kG\larr 0\,.
$$
 In particular, if $\,E\,$ is a cyclic configuration, then, using the
 notation of Section~\ref{sec2}, $\,\om_E\iso\kO_E\iso\kG(\oO,1,1)\,$
 and if $\,\kG=\kG(\fD,m,\la)\,$, there is an exact sequence
$$
 0\larr \kG\larr \kE \larr \kG\larr 0\,,
$$
 where
$$
 \kE\*_{\kO_E}\kG(\fD,m,\la)\iso
 \begin{cases}
 \kG(\fD,2,\la) &\text{if }\,m=1\,,\\
 \kG(\fD,m+1,\la)\+\kG(\fD,m-1,\la) &\text{if }\, m>1\,,
 \end{cases}
$$
 that is easily recognized as \AR{} (if $\,\chr\fK=0\,$ it follows
 immediately from \cite{ar1}, and the combinatoric description of
 $\VB(X)$ in \cite{dg2} does not depend on the characteristic). In
 particular, the \AQ{} of the category of vector bundles over a cyclic
 configuration is a disjoint union of \emph{homogeneous tubes}
$$
 \kT(\fD,\la):\quad \kG(\fD,1,\la) \toto \kG(\fD,2,\la) \toto
 \kG(\fD,3,\la) \toto \dots 
$$
 with $\,\tau\kG=\kG\,$ for every indecomposable vector bundle $\,\kG\,$.
 Here $\,\fD\,$ runs through aperiodic $\,s$-sequences and
 $\,\la\in\fK^*\,$.

 Passing from vector bundles to \Cm s, we use the following result of
 Kahn:
\begin{prop}
\label{kahn7}
 {\rm \cite[Theorem 3.1]{kah}}\ If $\,X\,$ is a minimally elliptic
 singularity, $\,Z\,$ its fundamental cycle and $\,0\to M\to N\to
 M\to0\,$ an \AR{} in $\,\CM(X)\,$, then the exact sequence $\,0\to
 R_ZM\to R_ZN\to R_ZM\to0\,$ is the direct sum of an \AR{}
 $\,0\to\kG\to\kF\to\kG\to0\,$ and a split sequence $\,0\to n\kO_Z\to
 2n\kO_Z\to n\kO_Z\to 0\,$.

\smallskip\noindent
\rm (Note that, as a minimally elliptic singularity is Gorenstein, any \AR{}
 in $\,\CM(X)\,$ is of the form $\,0\to M\to N\to M\to0\,$,
 cf. \cite[Theorem I.3.1]{ar1}.)
\end{prop}

\begin{corol}
\label{ars}
 All \AR s of \Cm s over a cusp singularity are the following:
\begin{align*}
& 0\to M(\fD,1,\la)\to M(\fD,2,\la)\to M(\fD,1,\la)\to 0\
 \textrm{ if }\, \fD\ne\fB\,\textrm{ or }\,\la\ne1\,;\\
& 0\to M(\fB,1,1)\to A\+M(\fB,2,1)\to M(\fB,1,1)\to 0\,;\\
& 0\to M(\fD,m,\la)\to M(\fD,m+1,\la)\+M(\fD,m-1,\la) \to
 M(\fD,m,\la) \to 0 \\& \textrm{ if }\, m>1\,.
\end{align*}
\end{corol}
 In particular, the \AQ{} of the category $\,\CM(X)\,$, where $\,X\,$ is a
 cusp singularity, is a disjoint union of homogeneous tubes
$$
 \dT(\fD,\la):\quad M(\fD,1,\la) \toto M(\fD,2,\la) \toto
 M(\fD,3,\la) \toto \dots 
$$
 where $\,\fD\ne\fB\,$ or $\,\la\ne1\,$, and one
 ``special'' tube 
$$
 \dT(\fB,1):\quad A\toto M(\fB,1,1) \toto M(\fB,2,1)
 \toto M(\fB,3,1) \toto \dots 
$$
 enlarged by the regular module $\,A\,$ that is both projective and
 ext-injective in this category (the latter means that $\,\ext^1_A(M,A)=0\,$
 for any \Cm\ $\,M\,$). Recall that in $\,M(\fD,m,\la)\,$
 always either $\,\fD>0\,$ or $\,\la\ne1\,$. Again $\,\tau M=M\,$ for
 all indecomposable modules $\,M\not\iso A\,$ (obviously, $\,\tau A\,$
 does not exist). 

\medskip
 Passing from the singularities $\,\rT_{pq2}\,$ to
 $\,\rT_{pq}\,$, one only has to note that Kn\"orrer's correspondence 
 maps an \AR{} $\,0\to M\to N\to M\to0\,$ either to an \AR{} $\,0\to
 \rest M\to \rest N\to\rest M\to 0 \,$ (if $\,M^\si\not\iso M\,$) or to
 the direct sum of \AR s $\,0\to\rest_2M\to N'\to\rest_1M\to0\,$ and
 $\,0\to \rest_1M\to N''\to\rest_2M\to0\,$ if $\,M\iso M^\si\,$
 (cf.~\cite[Corollary 2.10]{kno}). In particular, considering the \AR{}
 $\, 0\to M(\fB,1,1)\to A\+M(\fB,2,1)\to M(\fB,1,1)\to 0\,$, we see
 that $\,A'\,$   must be in the middle term of the \AR{} starting
 either from $\,N_1(\fB,1,1)\,$ or from $\,N_2(\fB,1,1)\,$. We choose
 the second possibility (it only depends on the notation). We also
 choose the notations in such a way that $\,N_i(\fD,m+1,\pm1)\,$ is
 a direct summand of the middle term of the \AR{} starting from
 $\,N_i(\fD,m,\pm1)\,$ ($\,i=1,2\,$).
 Then the \AQ{} is in this case a disjoint union of homogeneous tubes 
$$
 \tilde\dT(\fD,\la):\quad N(\fD,1,\la) \toto N(\fD,2,\la) \toto
 N(\fD,3,\la) \toto \dots 
$$
 for non-$\,\si$-symmetric $\,\fD\,$ or for $\,\la\ne1\,$; regular
 tubes $\,\tilde\dT(\fD,\pm1)\,$ of period $2$ for $\,\si$-symmetric
 $\,\fD\,$ ($\,\fD\ne\fB\,$ or $\,\la=-1\,$):
\nopagebreak\\
\setlength{\unitlength}{1cm}
\begin{picture}(15,3)
 \put(0,2){\parbox{2cm}{$N_1(\fD,1,\pm1)$}}
 \put(3.5,2){\parbox{2cm}{$N_1(\fD,2,\pm1)$}}
 \put(7,2){\parbox{2cm}{$N_1(\fD,3,\pm1)$}}
 \put(0,1){\parbox{2cm}{$N_2(\fD,1,\pm1)$}}
 \put(3.5,1){\parbox{2cm}{$N_2(\fD,2,\pm1)$}}
 \put(7,1){\parbox{2cm}{$N_2(\fD,3,\pm1)$}}
 \put(2.5,2.15){\vector(1,0){.7}} \put(6,2.15){\vector(1,0){.7}}
 \put(9.5,2.15){\vector(1,0){.7}} 
 \put(2.5,1.15){\vector(1,0){.7}} \put(6,1.15){\vector(1,0){.7}}
 \put(9.5,1.15){\vector(1,0){.7}} 
 \put(10.5,2){\parbox{1cm}{$\,\cdots\,$}}
 \put(10.5,1){\parbox{1cm}{$\,\cdots\,$}}
 \put(3.2,2){\vector(-1,-1){.72}}\put(3.2,1.3){\vector(-1,1){.72}}
 \put(6.7,2){\vector(-1,-1){.72}}\put(6.7,1.3){\vector(-1,1){.72}}
 \put(10.2,2){\vector(-1,-1){.72}}\put(10.2,1.3){\vector(-1,1){.72}}
\end{picture}\\
 and one ``special'' tube $\,\tilde\dT(\fB,1)\,$ of period $2$ enlarged
 by the regular module $\,A'\,$:\\ 
\setlength{\unitlength}{1cm}
\begin{picture}(15,3)
 \put(1,2){\parbox{1.6cm}{$N_1(\fB,1,1)$}}
 \put(4.3,2){\parbox{2cm}{$N_1(\fB,2,1)$}}
 \put(7.6,2){\parbox{2cm}{$N_1(\fB,3,1)$}}
 \put(1,1){\parbox{1.6cm}{$N_2(\fB,1,1)$}}
 \put(4.3,1){\parbox{2cm}{$N_2(\fB,2,1)$}}
 \put(7.6,1){\parbox{2cm}{$N_2(\fB,3,1)$}}
 \put(3.2,2.15){\vector(1,0){.7}} \put(6.5,2.15){\vector(1,0){.7}}
 \put(9.8,2.15){\vector(1,0){.7}} 
 \put(3.2,1.15){\vector(1,0){.7}} \put(6.5,1.15){\vector(1,0){.7}}
 \put(9.8,1.15){\vector(1,0){.7}} 
 \put(10.9,2){\parbox{1cm}{$\,\cdots\,$}}
 \put(10.9,1){\parbox{1cm}{$\,\cdots\,$}}
 \put(3.9,2){\vector(-1,-1){.72}}\put(3.9,1.3){\vector(-1,1){.72}}
 \put(7.2,2){\vector(-1,-1){.72}}\put(7.2,1.3){\vector(-1,1){.72}}
 \put(10.5,2){\vector(-1,-1){.72}}\put(10.5,1.3){\vector(-1,1){.72}}
 \put(0,1.5){\parbox{.3cm}{$A'$}}
 \put(.9,1.2){\vector(-4,3){.4}}\put(.5,1.8){\vector(4,3){.4}}
\end{picture}
 The action of the \AT{} $\,\tau\,$ is trivial in homogeneous tubes and
 coincides with the obvious axis symmetry in those of period $2$.
 
 The procedure is quite the same when passing to hypersurface
 singularities $\,\rT_{pqr}\,$ of higher dimensions. Namely, again the
 mapping $\,\syz\,$ transforms an \AR{} $\,0\to M\to N\to M\to 0\,$
 either to an \AR{} $\,0\to\syz M\to\syz N\to \syz M\to0\,$ (if
 $\,M\not\iso\Om M\,$) or to the direct sum of \AR s $\,0\to\syz_2M\to
 N'\to\syz_1M\to0\,$ and $\,0\to\syz_1M\to N''\to\syz_2M\to0\,$. 

 For a hypersurface singularity of type $\,\rT_{pqr}\,$
 ($\,1/p+1/q+1/r<1\,$), one can use \emph{Kn\"orrer's periodicity}
 \cite{kno,yos} to describe the part of the \AQ{} not containing the free
 module and the result of S{\o}lberg \cite{sol} on the position of the
 free module in an \AR. Then one obtains that if the dimension
 $\,n=k+3\,$ is even, the \AR{} is the same as for the corresponding surface
 singularity, with the only exception that the tube containing the
 free module begins not from two arrows but from $\,2^{n/2}\,$ ones
 (half of them starting and half of them ending at $\,A\,$). If
 $\,n\,$ is odd, the \AR{} looks like that for the curve case, except of the
 tube containing the free module, where to each arrow starting or ending
 at $\,A\,$ should be added $\,2^{[\frac{n-1}2]}\,$ ones such
 that each next arrow goes to the opposite direction with respect to the
 preceding one.

\section{Summary and conjectures}

 We can now summarize the known results on \Com{} types of isolated
 \Com\ singularities. This is done in Table 1. To make it more uniform, we call a curve
 singularity of type $\,\rT_{pq}\,$ \emph{simple elliptic} if
 $\,1/p+1/q=1/2\,$ and a \emph{cusp} if $\,1/p+1/q<1/2\,$. In the same
 way, we call a hypersurface singularity of type $\,\rT_{pqr}\,$
 \emph{simple elliptic} if $\,1/p+1/q+1/r=1\,$ and a \emph{cusp} if
 $\,1/p+1/q+1/r<1\,$, though it seems not to be the usual practice. 

 Moreover, there is some evidence that these results are
 \emph{complete}, that is, all the remaining singularities are \Com{}
 wild. We formulate the corresponding conjectures as well as one
 related to non-isolated singularities. Here $\,(X,x)\,$ denotes
 \emph{any} \Com\ singularity of an algebraic variety over a field of
 characteristic $\,0\,$ and \Com{} type refers to its complete local ring
 $\,A\,$. 

\begin{table}
 \begin{center}
 Table 1.\\\medskip \sf Cohen--Macaulay types of singularities
 \end{center}
\bigskip 
\begin{center}
\begin{tabular}{|c||c|c|c|}
\hline &&&\\
 CM type & curves & surfaces & hypersurfaces \\
&&&\\ \hline\hline &&&\\
 finite & dominate  & quotient & simple \\
 & simple &&\\
&&&\\ \hline &&&\\
 tame & dominate  & simple & simple \\
 bounded & simple    & elliptic- & elliptic \\
      & elliptic & quotient & (no other\,?) \\
 && (no other\,?) & \\
&&&\\ \hline &&&\\
  tame & dominate  & cusp- & cusp \\
 unbounded & cusp    & quotient & (no other\,?) \\
      && (no other\,?) & \\
&&&\\ \hline &&&\\
 wild & all other & all other\,? & all other\,? \\
&&&\\ \hline 
\end{tabular}
\end{center}
\end{table}

\begin{conjecture}
 In the following cases the ring $\,A=\hat\kO_{X.x}\,$ is \Com{} wild:
\begin{enumerate}
 \item
  $\,(X,x)\,$ is a surface singularity which is neither a quotient nor an
  elliptic-quotient.
 \item
  $\,(X,x)\,$ is a hypersurface singularity, which is neither a simple
  one nor of type $\,\rT_{pqr}\,$.
 \item
  $\,(X,x)\,$ is a non-isolated singularity with the dimension of the
  singular locus greater than $\,1\,$.
\end{enumerate}
\end{conjecture}
 If this conjecture is true, it yields a complete description of
 isolated \Com{} tame surface and hypersurface singularities together
 with a classification of their indecomposable \Cm s. 

\medskip
\begin{remk}\label{expgrow}
 All known examples of \Com\ tame unbounded singularities, in particular those from Table 1,
 are actually of exponential growth.  It seems very plausible that it is always so. Nevertheless,
 just as in the case of finite dimensional algebras, it can only be shown \emph{a posteriori},
 when one has a description of modules. We do not see any ``natural'' way to prove this conjecture
 without such calculations. 
\end{remk}

\medskip
\begin{remk}
\label{anal}
 In the complex analytic case, Artin's Approximation Theorem \cite{art} implies
 that the list of \Cm s (Corollary~\ref{cusp1}) remains the same if
 $\,A\,$ denotes the ring of germs of analytic functions on a cusp
 singularity. The lifting of families (Propositions~\ref{kahn5} and
 \ref{kahn6}) is more cumbersome. We do only claim that, for each
 point $\,t\in T\,$, a lifting is possible over a neighbourhood $\,U\,$ of
 $\,t\,$ in $\,T\,$. Combined with the uniqueness assertion, just as in Remark~\ref{glob},
 it gives a lifting of an appropriate family to the
 universal covering $\,\tilde T\,$ of $\,T\,$. If $\,T\,$ is a smooth
 curve, so is $\,\tilde T\,$, therefore the results on tameness from
 Section~\ref{sec4} remain valid. On the other hand, in the case of cusps it seems credible
 that the families $\,\kG(\fD)\,$ from the proof of Corollary~\ref{cusp2}
 can actually be lifted over $\,T\,$, just as in \cite{kah} for simple elliptic
 case. Moreover, in the case (iv) of Remark~\ref{re1}, where the base is a finite dimensional algebra,
 Artin's Approximation Theorem can also be applied, so
 Theorem~\ref{minel} remains valid in the analytic case too.
\end{remk}


\begin{thebibliography}{55}

\bibitem{arn}
 V.~I.~Arnold, S.~M.~Gusein-Zade and A.~N.~Varchenko.
 \emph{Singularities of Differentiable Maps.} Vol.\,I. Birkh\"auser,
 Boston, 1985.
 
\bibitem{art}  
 M.~Artin. Algebraic approximation of structures over complete local
 rings. Publ.~math.~IHES, 36 (1969), 23--58.

\bibitem{au2}
 M.~Auslander. On the purity of the branch loci. Amer.~J.~Math. 84
(1962), 116--125.
 
\bibitem{aus} 
 M.~Auslander. Rational singularities and almost split sequences.
 Trans. Amer. Math. Soc. 293 (1986), 511-531

\bibitem{ab}
 M.\,Auslander~and~D.\,A.\,Buchsbaum. On ramification theory in noetherian rings.
 Amer.~J.~Math. 81 (1959), 749--765.

\bibitem{ar1}
 M.~Auslander and I.~Reiten.
 Almost split sequences in dimension two. Adv. Math. 66 (1987),
 88--118. 

\bibitem{dr1} 
 Y.\,A.~Drozd. Tame and wild matrix problems. In: \emph{Representations
 and Quadratic Forms}, Institute of mathematics, Kiev, 1979,
39--74. (English translation in: Amer. Math. Soc. Translations, 128
(1986), 31--55.)

\bibitem{dg1} 
 Y.\,A.~Drozd~and~G.-M.Greuel. Cohen--Macaulay module type. Compositio Math.
89 (1993), 315--338.

\bibitem{dg2}  
 Y.\,A.~Drozd~and~G.-M.Greuel.  Tame and wild projective
curves and classification of vector bundles. J.~Algebra, 246 (2001) 1--54.

\bibitem{dk}
 Y.\,A.~Drozd~and~V.\,V.~Kirichenko. \emph{Finite Dimensional Algebras.}
 Springer, New York, 1994.

\bibitem{drr} 
 Y.\,A.~Drozd and A.\,V.~Roiter.
 Commutative  rings  with  a finite number of integral indecomposable
 representations, Izvestia Acad. Sci. USSR, 31 (1967) 783--798. 

\bibitem{esn} 
 H.~Esnault. Reflexive modules on quotient surface
  singularities. J. Reine Angew, Math. 362 (1985), 63--71. 

\bibitem{grk}
 G.-M.~Greuel and H.~Kn\"orrer.
 Einfache Kurvensingularit\"aten und torsionfreie
 Moduln. Math. Ann. 270 (1985), 417--425.

\bibitem{gr} A.~Grothendieck. \emph{\'El\'ements de G\'eom\'etrie Alg\'ebrique,\,IV}. 
 Publ. Math. I.H.E.S. 24 (1965).

\bibitem{har} 
 R.~Hartshorne. \emph{Algebraic Geometry}. Springer, New York, 1977.

\bibitem{hir} 
 F.~Hirzebruch. The Hilbert modular group, resolution of 
 the singularities at the cusps and related problems.
 Seminaire Bourbaki~1970/71, Expos\'e~396.
  
\bibitem{jac}
 H.~Jacobinski.
 Sur les ordres commutatifs avec un nombre fini de r\'eseaux
 indecomposables. Acta Math. 118 (1976), 1--31.
  
\bibitem{kah} 
 C.~Kahn. Reflexive modules on minimally elliptic singularities.
 Math.Ann. 285 (1989), 141-160.

\bibitem{kar}
 U.~Karras. 
 Dissertation, Universit\"at Bonn, 1973.

\bibitem{kaw} 
 Y.~Kawamata. Crepant blowing-up of 3-dimensional canonical
 singularities and its application to degenerations of surfaces.
 Ann. Math. 127 (1988), 93--163.

\bibitem{kno} 
 H.~Kn\"orrer. Maximal \Cm s on hypersurface singularities,
 I. Invent. Math. 88 (1987), 153--164.

\bibitem{lau}
 H.~Laufer.
 On minimally elliptic singularities.
 Amer.~J.~Math. 99 (1975), 1257--1295.
 
\bibitem{lip}
 J.~Lipman. Rational singularities, with applications to
 algebraic surfaces and unique factorization. Publ.~Math.~IHES 36
 (1969), 195--279.

\bibitem{ser} 
 J.-P.~Serre. \emph{Alg\`ebre locale, multiplicit\'es}.
 Lecture Notes in Math. 11, Springer, Berlin, 1975.

\bibitem{sol}
 \O.~Solberg. Hypersurface singularities of finite \Com{}
 type. Proc. London. Math. Soc. 58 (1989) 258--280.

\bibitem{yos}
 Y.\,Yoshino.
 \emph{Cohen--Macaulay Modules over Cohen--Macaulay Rings}. Cambridge
 University Press, Cambridge, 1990.
 
\end{thebibliography}
\end{document}